\documentclass[aop,citesort,MSNbibl,dvips]{arximspdf}
\usepackage{mathrsfs}
\usepackage{dcolumn}
\usepackage{graphicx}

\doi{10.1214/09-AOP489}
\volume{38}
\issue{2}
\pubyear{2010}
\firstpage{570}
\lastpage{604}

\makeatletter

\newcolumntype{d}[1]{D{.}{.}{#1}}

\newtheorem{theorem}{Theorem}[section]
\newtheorem{proposition}[theorem]{Proposition}
\newtheorem{lemma}[theorem]{Lemma}

\newproclaim{definition}[theorem]{Definition}
\newproclaim{remark}[theorem]{Remark}
\newproclaim{example}[theorem]{Example}

\newproclaim{Case}{Case}

\makeatother

\begin{document}
\begin{frontmatter}

\title{Coverage processes on spheres and
condition numbers for linear programming}
\runtitle{Coverage processes and condition numbers}

\begin{aug}
\author[A]{\fnms{Peter} \snm{B\"{u}rgisser}\corref{}\thanksref{t1}\ead[label=e1]{pbuerg@upb.de}},
\author[B]{\fnms{Felipe} \snm{Cucker}\thanksref{t2}\ead[label=e2]{macucker@cityu.edu.hk}} and
\author[C]{\fnms{Martin} \snm{Lotz}\ead[label=e3]{lotz@maths.ox.ac.uk}}
\runauthor{P. B\"{u}rgisser, F. Cucker and M. Lotz}
\affiliation{University of Paderborn, City University of Hong Kong and
Oxford University}
\address[A]{P. B\"{u}rgisser\\
Institute of Mathematics\\
University of Paderborn\\
33098 Paderborn\\
Germany\\
\printead{e1}} 
\address[B]{F. Cucker\\
Department of Mathematics\\
City University of Hong Kong\\
Kowloon Tong\\
Hong Kong\\
\printead{e2}}
\address[C]{M. Lotz\\
Mathematical Institute\\
University of Oxford\\
24-29 St. Giles'\\
Oxford OX1 3LB\\
England\\
\printead{e3}}
\end{aug}

\thankstext{t1}{Supported in part by DFG Grant BU 1371/2-1.}

\thankstext{t2}{Supported in part by CityU GRF Grant CityU 100808.}

\pdfauthor{Peter Burgisser, Felipe Cucker, Martin Lotz}

\received{\smonth{4} \syear{2008}}
\revised{\smonth{3} \syear{2009}}

%
\begin{abstract}
This paper has two agendas. Firstly, we exhibit new results
for coverage processes.
Let $p(n,m,\alpha)$ be the probability that $n$ spherical
caps of angular radius $\alpha$ in $S^m$ do not cover
the whole sphere $S^m$.
We give an exact formula for $p(n,m,\alpha)$
in the case $\alpha\in[\pi/2,\pi]$
and an upper bound for $p(n,m,\alpha)$ in the case $\alpha\in
[0,\pi/2]$
which tends to $p(n,m,\pi/2)$ when $\alpha\to\pi/2$.
In the case $\alpha\in[0,\pi/2]$
this yields upper bounds for the expected number of spherical caps
of radius $\alpha$ that are needed to cover $S^m$.

Secondly, we study the condition number ${\mathscr C}(A)$ of the
linear programming feasibility problem
$\exists x\in\mathbb R^{m+1} Ax\le0, x\ne0$
where $A\in\mathbb R^{n\times(m+1)}$
is randomly chosen according to the standard normal distribution.
We exactly determine the distribution of ${\mathscr C}(A)$ conditioned
to $A$ being feasible and
provide an upper bound on the distribution function
in the infeasible case. Using these results,
we show that
$\mathbf{E}(\ln{\mathscr C}(A))\le2\ln(m+1) + 3.31$
for all $n>m$, the sharpest bound for this
expectancy as of today.
Both agendas are related through a result
which translates between coverage and condition.
\end{abstract}

%

\end{frontmatter}

\section{Introduction}\label{sec1}

\subsection{Coverage processes on spheres}

\begin{quote}
\textit{One of the oldest problems in the theory of coverage processes
is that
of calculating the chance that a given region is completely covered by a
sequence of random sets. Unfortunately there is only a small number of
useful circumstances where this probability may be calculated explicitly.}
(Hall \cite{Hall88}, Section 1.11.)
\end{quote}

In 1897 Whitworth \cite{Whit1897} considered the following problem.
Assume we place $n$ arcs of angular radius $\alpha$ in the unit
circle $S^1$,
whose centers are independently and randomly chosen from the
uniform distribution in $S^1$. What is the probability that these
arcs do not cover $S^1$?

Whitworth's problem is arguably at the origin of the theory of
coverage processes. It was not until 1939 that an answer to the
problem was given when Stevens \cite{Stev39} showed that the probability
in question is
%
%
\begin{equation}\label{eq:stevens}
\sum_{j=1}^k (-1)^{j+1}
\pmatrix{n\cr j}\biggl(1-\frac{j\alpha}{\pi} \biggr)^{n-1},
\end{equation}
where $k=\lfloor\frac{\pi}{\alpha}\rfloor$. Extensions of this
result to
other quantities related with random arcs in $S^1$ are given
in \cite{Siegel79}. Extensions to random arcs with different lengths
are given in \cite{Dvor56,Kahane59} and in \cite{SieHo82}
where an exact formula for the probability above is given for randomly
placed arcs having random independent size.

The extension of the original problem in $S^1$
to the two-dimensional unit sphere $S^2$ was considered by
Moran and Fazekas de St. Groth \cite{MoFa62}. Let $p(n,\alpha)$
denote the probability that $n$ spherical caps of angular radius
$\alpha$,
and centers randomly and independently
chosen from the uniform distribution on $S^2$ do not cover
$S^2$. Moran and Fazekas de St. Groth exhibited an approximation
of $p(n,\alpha)$, and numerically estimated this quantity for
$\alpha=53^\circ26'$ (a value arising in a biological problem motivating
their research). Shortly thereafter, Gilbert \cite{Gilbert65}
showed the bounds
%
%
\begin{equation}\label{eq:Gilbert}
(1-\lambda)^n \leq p(n,\alpha)
\leq\tfrac43 n(n-1)\lambda(1-\lambda)^{n-1},
\end{equation}
where $\lambda= (\sin\frac\alpha{2} )^2=\frac12(1-\cos\alpha)$
is the fraction of the surface
of the sphere covered by each cap. In addition, Gilbert conjectured
that, for $n\to\infty$, $p(n,\alpha)$ satisfies the asymptotic equivalence
\[
p(n,\alpha)\approx n(n-1)\lambda^2
(1-\lambda^2 )^{n-1}.
\]
This conjecture was proven by Miles \cite{Miles69} who also found
an explicit expression (cf. \cite{Miles68a})
for $p(n,\alpha)$ if $\alpha\in[\pi/2,\pi]$,
namely
%
%
\begin{eqnarray}\label{eq:Miles}
p(n,\alpha) &=&\pmatrix{n\cr2}\int_0^{\pi-\alpha}
\sin^{2(n-2)}(\theta/2)
\sin(2\theta) \,d\theta\nonumber\\[-8pt]\\[-8pt]
&&{} + \frac34
\pmatrix{n\cr3}\int_0^{\pi-\alpha}
\sin^{2(n-3)}(\theta/2)
\sin^3\theta \,d\theta.\nonumber
\end{eqnarray}
More on the coverage problem for $S^1$ and $S^2$ can be
found in \cite{solo78}.
Extensions of these results to the unit sphere $S^m$
in $\mathbb R^{m+1}$ for $m>2$ are scarce. Let $p(n,m,\alpha)$ be the
probability that $n$ spherical
caps of angular radius $\alpha$ in $S^m$ do not cover $S^m$. That is,
for $\alpha\in[0,\pi]$, and $a_1,\ldots,a_n$ randomly and independently
chosen points in $S^m$ from the uniform distribution, define
\[
p(n,m,\alpha):=
\mathsf{Prob}\Biggl\{S^m\neq\bigcup_{i=1}^n\operatorname{\mathsf{cap}}(a_i,\alpha)
\Biggr\},
\]
where $\mathsf{cap}(a,\alpha)$ denotes the spherical cap of angular radius
$\alpha$
around $a$. It can easily be seen that for $n\leq m+1$ and $\alpha\leq
\pi/2$
we have $p(n,m,\alpha)=1$.
Moreover, Wendel \cite{Wendel62} has shown that
%
%
\begin{equation}\label{eq:wendel}
p(n,m,\pi/2)=2^{1-n}\sum_{k=0}^{m}\pmatrix{n-1\cr k}.
\end{equation}
Furthermore, a result by Janson \cite{Janson86} gives an asymptotic
estimate of $p(n,m,\alpha)$ for $\alpha\to0$. Actually, Janson's article
covers a situation much more general than fixed radius caps on a sphere
and it was preceded by a paper by Hall \cite{Hall85} where bounds for
the coverage probability were shown for the case of random spheres on
a torus.

A goal of this paper is to extend some of the known results
for $S^1$ and $S^2$ to higher dimensions. To describe our results
we first introduce some notation. We denote by
\[
\mathcal O_m := \mathsf{vol}_{m}(S^m)
= \frac{2\pi^{({m+1})/{2}}}{\Gamma(({m+1})/{2} )}
\]
the $m$-dimensional volume of the sphere $S^m$. Also,
for $t\in[0,1]$, denote the relative volume of a cap
of radius $\arccos t\in[0,\pi/2]$ in $S^m$ by $\lambda_{m}(t)$.
It is well known that
%
%
\begin{equation}\label{eq:volcap}
\lambda_{m}(t) = \frac{\mathcal O_{m-1}}{\mathcal O_{m}}
\int_0^{\arccos t} (\sin\theta)^{m-1} \,d\theta.
\end{equation}

Our results are formulated in terms of a family of numbers $C(m,k)$
defined for $1\leq k\leq m$.
These numbers are defined in Section \ref{se:nJac} and
studied in Section \ref{se:C(m,k)}
where we give bounds on $C(m,k)$ and derive a closed form for
$k\in\{1,m-1,m\}$.
Furthermore, we will show that, for each $m$, the $C(m,k)$
can be obtained as the solution of a system
of linear equations which easily allows us to
produce a table for their values
(cf. Table \ref{table:1}).

A main result in this paper is the following.
\begin{theorem}\label{th:main_coverage}
Let $n >m\ge1$,
$\alpha\in[0,\pi]$, and $\varepsilon=\cos(\pi-\alpha)$.
For $\alpha\in[\frac{\pi}{2},\pi]$
\[
p(n,m,\alpha)
= \sum_{k=1}^m \pmatrix{n\cr{k+1}} C(m,k)
\int_{\varepsilon}^{1} t^{m-k}(1-t^2)^{km/2-1}
\lambda_{m}^{n-k-1}(t) \,dt
\]
and for $\alpha\in[0,\frac{\pi}{2})$ we have
\begin{eqnarray*}
p(n,m,\alpha)&\leq&
\frac{\sum_{k=0}^{m}{n-1\choose k}}
{2^{n-1}}\\
&&{} + \pmatrix{n\cr{m+1}} C(m,m)
\int_0^{|\varepsilon|} (1-t^2)^{({m^2-2})/{2}}
\bigl(1-\lambda_{m}(t) \bigr)^{n-m-1} \,dt.
\end{eqnarray*}
\end{theorem}

We remark that this formula, for $\alpha\in[\pi/2,\pi]$ and $m=2$,
is identical to the one given by Miles (\ref{eq:Miles}). Also,
for $\alpha\in[0,\pi/2]$ and $m=1$, our upper bound for
$p(n,1,\alpha)$ coincides with the first term in Steven's
formula (\ref{eq:stevens}) (cf. Remark \ref{rem:stevens} below).

%
%
\begin{table}
\tablewidth=250pt
\caption{A few values for $C(m,k)$}\label{table:1}
\begin{tabular*}{\tablewidth}{@{\extracolsep{\fill}}lccccd{2.4}c@{}}
\hline
$\bolds{k\setminus m}$ & \textbf{1} & \textbf{2} & \textbf{3} & \textbf{4} & \textbf{5} & \textbf{6}\\
\hline
1 & $\frac2{\pi}$ & 2 & 5.0930 & 12 & 27.1639 &60\\
2 && $3/4$ & 3.9317 & $477/32$ & 49.5841 & $78795/512$\\
3 &&& 0.6366 & $39/8$ & 25.1644 & $897345/8192$ \\
4 &&&& $15/32$ & 4.8525 & $132225/4096$ \\
5 &&&&& 0.3183 & $4335/1024$ \\
6 &&&&&& $105/512$\\
\hline
\end{tabular*}
\end{table}

We may use Theorem \ref{th:main_coverage}, together with the bound on
the $C(m,k)$, to derive bounds for the expected value of $N(m,\alpha
)$, the
number of random caps of radius $\alpha$ needed to cover $S^m$. The
asymptotic behavior of $N(m,\alpha)$ for $\alpha\to0$ has been
studied by
Janson \cite{Janson86}. Otherwise, we have not found any bound for
$\mathbf{E}(N(m,\alpha))$ in the literature.
\begin{theorem}\label{cor:2}
For $\alpha\in(0,\frac{\pi}{2}]$ we have
\[
\mathbf{E}(N(m,\alpha))\leq3m + 2 + \sqrt{m} (m+1) \cos(\alpha)
\lambda_m(\cos (\alpha))^{-2}
\biggl(\frac{1}{2\lambda_m(\cos(\alpha))} \biggr)^m.
\]
\end{theorem}

\subsection{Polyhedral conic systems and their condition}

Among the number of interrelated problems collectively known
as linear programming, we consider the following two.

\subsubsection*{Feasibility of polyhedral conic systems
(FPCS)}
Given a matrix $A\in\mathbb R^{n\times(m+1)}$, decide whether there
exists a nonzero $x\in\mathbb R^{m+1}$ such that $Ax\leq0$ (componentwise).

\subsubsection*{Computation of points in polyhedral cones
(CPPC)}
Given a matrix $A\in\mathbb R^{n\times(m+1)}$ such that
${\mathscr S}=\{x\in\mathbb R^{m+1}\mid Ax<0\}\neq\varnothing$, find
$x\in{\mathscr S}$.\vspace*{12pt}

By scaling we may assume without loss of generality that the rows
$a_1,\ldots,a_n$ of $A$ have Euclidean norm one and interpret
the matrix $A$ as a point in $(S^m)^n$.
We say that the elements of the set
%
%
\begin{equation}\label{eq:F}
\mathcal{F}_{n,m} :=\{A\in(S^m)^n \mid
\exists x\in S^m \langle a_1,x\rangle\le0,\ldots,\langle a_n,x\rangle
\le0\},
\end{equation}
are \textit{feasible}. Similarly, we say that the elements in
%
%
\begin{equation}\label{eq:Fcirc}
\mathcal{F}_{n,m}^\circ:=\{A\in(S^m)^n \mid
\exists x\in S^m \langle a_1,x\rangle<0,\ldots,\langle a_n,x\rangle
<0\}
\end{equation}
are \textit{strictly feasible}. Elements in $(S^m)^n\setminus\mathcal{F}_{n,m}$
are called \textit{infeasible}. Finally, we call \textit{ill-posed} the
elements in
$\Sigma_{n,m}:=\mathcal{F}_{n,m}\setminus\mathcal{F}_{n,m}^{\circ}$.

For several iterative algorithms solving the two problems
above, it has been observed that the number of iterations required
by an instance $A$ increases with the quantity
\[
{\mathscr C}(A)=\frac{1}{\mathsf{dist}(A,\Sigma_{n,m})},
\]
(here $\mathsf{dist}$ is the distance with respect to an appropriate metric;
for the precise definition we refer to Section \ref{se:prel-poly}).
This quantity, known as the \textit{GCC-condition number} of
$A$ \cite{Goff80,ChC00}, occurs together with the dimensions $n$ and $m$
in the theoretical analysis (for both complexity and accuracy) of
the algorithms mentioned above.
For example, a primal-dual interior-point method is
used in \cite{CP01} to solve (FPCS) within
%
%
\begin{equation}\label{eq:cotaCP}
\mathcal O\bigl(\sqrt{m+n}\bigl(\ln(m+n)+\ln{\mathscr C}(A)\bigr) \bigr)
\end{equation}
iterations. The Agmon--Motzkin--Sch\"{o}nberg relaxation
method\footnote{This method gives also the context in which
${\mathscr C}(A)$ was first studied, although in the feasible
case only \cite{Goff80}.} \cite{Agmon,Motzkin} or the perceptron
method \cite{Rosenblatt} solve (CPPC) in a number of iterations
of order $\mathcal O({\mathscr C}(A)^2)$ (see Appendix B of \cite
{ChCH05} for a
brief description of this).

The complexity bounds above, however, are of limited use since, unlike
$n$ and~$m$, ${\mathscr C}(A)$ cannot be directly read from $A$. A way to
remove ${\mathscr C}(A)$ from these bounds consists in trading worst-case
by average-case analysis. To this end, one endows the space $(S^m)^n$
of matrices $A$ with a probability measure and studies ${\mathscr C}(A)$
as a random variable with the induced distribution.
In most of these works, this measure is the unform one in $(S^m)^n$
(i.e., matrices $A$ are assumed to have its $n$ rows independently
drawn from the uniform distribution in $S^m$).

Once a measure has been set on the space of matrices [and in what follows
we will assume the uniform measure
in $(S^m)^n$], an estimate on $\mathbf{E}(\ln{\mathscr C}(A))$ yields
bounds on the average
complexity for (FPCS) directly from (\ref{eq:cotaCP}). For
(CPPC) the situation is different since it is
known \cite{ChCH05}, Corollary 9.4, that $\mathbf{E}({\mathscr
C}(A)^2)=\infty$. Yet,
an estimate for $\varepsilon>0$ on
\[
\mathsf{Prob}\{{\mathscr C}(A)\geq1/\varepsilon\mid A\in
\mathcal{F}_{n,m}\}
\]
yields bounds on the probability that the relaxation or perceptron
algorithms will need more than a given number of iterations.
Efforts have therefore been devoted to compute the expected value (or the
distribution function) of ${\mathscr C}(A)$ for random matrices $A$.

Existing results for these efforts are
easily summarized.
A bound for $\mathbf{E}(\ln{\mathscr C}(A))$ of the form $\mathcal
O(\min\{n, m\ln n\})$
was shown in \cite{ChC01}. This bound was improved \cite{CW01} to
$\max\{\ln m, \ln\ln n \}+\mathcal O(1)$ assuming that $n$ is moderately
larger than $m$. Still, in \cite{ChCH05}, the asymptotic behavior
of both ${\mathscr C}(A)$ and $\ln{\mathscr C}(A)$ was exhaustively
studied, and these results were extended in \cite{HM06} to
matrices $A\in(S^m)^n$ drawn from distributions
more general than the uniform.
Independently of this stream of results, in \cite{DST}, a
smoothed analysis for a related condition number is performed
from which it follows that $\mathbf{E}(\ln{\mathscr C}(A))=\mathcal
O(\ln n)$.

Our second set of results adds to the line
of research above. First, we
provide the exact distribution of ${\mathscr C}(A)$ conditioned to
$A$ being feasible and a bound on this distribution for the
infeasible case.
\begin{theorem}\label{th:main_lp}
Let $A\in(S^m)^n$ be randomly chosen
from the uniform distribution in $(S^m)^n$, $n>m$.
Then, for $\varepsilon\in(0,1]$, we have
\begin{eqnarray*}
&&\mathsf{Prob}\{{\mathscr C}(A)\geq1/\varepsilon\mid A\in
\mathcal{F}_{n,m} \}\\
&&\qquad= \frac{2^{n-1}}{\sum_{k=0}^{m}{{n-1}\choose k}}
\sum_{k=1}^m \pmatrix{n\cr{k+1}} C(m,k)
\\
&&\qquad\quad\hspace*{52.65pt}{}\times
\int_0^{\varepsilon} t^{m-k}(1-t^2)^{km/2-1}
\lambda_{m}(t)^{n-k-1} \,dt,
\\
&&\mathsf{Prob}\{{\mathscr C}(A)\geq1/\varepsilon\mid A\notin
\mathcal{F}_{n,m} \}\\
&&\qquad\leq\frac{2^{n-1}}{\sum_{k=m+1}^{n-1}{{n-1}\choose k}}
\pmatrix{n\cr{m+1}} C(m,m)\\
&&\qquad\quad{}\times
\int_0^{\varepsilon} (1-t^2)^{({m^2-2})/{2}}
\bigl(1-\lambda_m(t) \bigr)^{n-m-1} \,dt.
\end{eqnarray*}
\end{theorem}

Second, we prove an upper bound on $\mathbf{E}(\ln{\mathscr C}(A))$ that
depends only on $m$, in sharp contrast with all the previous
bounds for this expected value.
\begin{theorem}\label{corol:1}
For matrices $A$ randomly chosen
from the uniform distribution in $(S^m)^n$ with $n>m$, we have
$\mathbf{E}(\ln{\mathscr C}(A))\leq2\ln(m+1) + 3.31$.
\end{theorem}

Note that the best previously established
upper bound for $\mathbf{E}(\ln{\mathscr C}(A))$ (for arbitrary
values of $n$ and $m$)
was $\mathcal O(\ln n)$. The bound $2\ln(m+1) + 3.31$ is not only sharper
(in that it is independent of $n$) but also more precise
(in that it does not rely on the $\mathcal O$ notation).\footnote
{Recently a
different derivation of a $\mathcal O(\ln m)$ bound for $\mathbf
{E}(\ln{\mathscr C}(A))$
was given in \cite{AmBu08}. However, this derivation does not provide
explicit estimates for the constant in the $\mathcal O$ notation.}

\subsection{Coverage processes versus condition numbers}

Theorems \ref{th:main_coverage} and \ref{th:main_lp} are not
unrelated. Our next result, which will be the first one we
will prove, shows a precise link between coverage processes
and condition for polyhedral conic systems.
\begin{proposition}\label{prop:equiv}
Let $a_1,\ldots,a_n$ be randomly chosen from the uniform
distribution in $S^m$.
Denote by $A$ the matrix with rows $a_1,\ldots,a_n$. Then,
setting $\varepsilon:=|{\cos}(\alpha)|$ for $\alpha\in[0,\pi]$,
we have
\[
p(n,m,\alpha)=\cases{
\mathsf{Prob}\biggl\{A\in\mathcal{F}_{n,m} \mbox{ and }
{\mathscr C}(A)\leq\dfrac1{\varepsilon} \biggr\},
&\quad if $\alpha\in[\pi/2,\pi]$,\cr
\dfrac{\sum_{k=0}^{m}{n-1\choose k}}
{2^{n-1}}\cr
\qquad{}+\mathsf{Prob}\biggl\{A\notin\mathcal{F}_{n,m}
\mbox{ and } {\mathscr C}(A)\geq\dfrac1{\varepsilon} \biggr\},
&\quad if $\alpha\in[0,\pi/2]$.}
\]
In particular,
$p(n,m,\pi/2)=\mathsf{Prob}\{A\in\mathcal{F}_{n,m}\}
=2^{1-n}\sum_{k=0}^{m}\pmatrix{n-1\cr k}$.
\end{proposition}

While Proposition \ref{prop:equiv} provides a dictionary between
the coverage problem in the sphere and the condition of polyhedral
conic systems, it should be noted that, traditionally, these problems
have not been dealt with together. Interest on the second focused
on the case of ${\mathscr C}(A)$ being large or, equivalently, on
$\alpha$ being close to $\pi/2$. In contrast, research on the first
mostly focused on asymptotics for either
small $\alpha$ or large $n$ (an exception being \cite{Wendel62}).

\section{Main ideas}\label{se:mainideas}

In this section we describe in broad strokes how the results presented
in the \hyperref[sec1]{Introduction} are arrived at. In a first step
in Section \ref{se:prel-poly}, we give a characterization of the GCC
condition number which establishes a link to covering problems thus
leading to a proof of Proposition \ref{prop:equiv}. We then proceed
by explaining the main ideas behind the proof of
Theorem \ref{th:main_lp}.

In all that follows, we will write $[n]=\{1,\ldots,n\}$ for $n\in
\mathbb N$.

\subsection{The GCC condition number and spherical caps}
\label{se:prel-poly}

\begin{figure}

\includegraphics{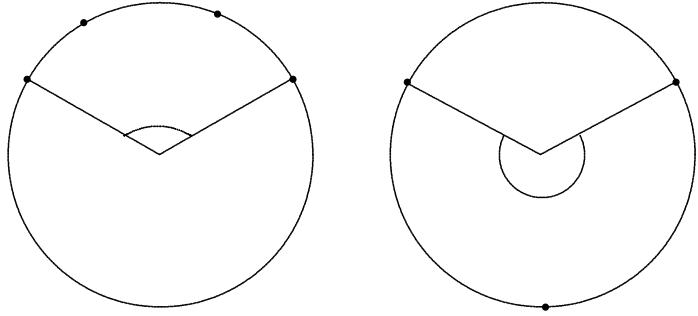}

\caption{A SIC with $\alpha\in(0,\pi/2)$ (left) and with $\alpha
\in(\pi
/2,\pi)$ (right).}
\label{fig:1}
\end{figure}

A key ingredient in what follows is a way of characterizing
the GCC condition number in terms of spherical caps.
For $p\in S^m$ and $\alpha\in[0,\pi]$, recall that
\[
\mathsf{cap}(p,\alpha) := \{x\in S^m\mid\langle p,x\rangle\ge\cos
\alpha\} .
\]
A \textit{smallest including cap} ($\mathrm{SIC}$)
for $A=(a_1,\ldots,a_n)\in(S^m)^n$ is a spherical cap
of minimal radius containing the points
$a_1,\ldots,a_n$.
If $p$ denotes its center, then its \textit{blocking set} is defined as
$\{i\in[n]\mid\langle p,a_i\rangle= \cos\alpha\}$
which can be seen as a set of ``active'' rows
(cf. Figure \ref{fig:1}).

%
%

A \textit{largest excluding cap} (LEC) for $A$ is the complement
of a smallest including cap for $A$.
Note that (by a compactness argument) a $\mathrm{SIC}$ always exists,
and while
there may be
several $\mathrm{SIC}$ for $A$, its radius is uniquely determined.
For the rest of this article, we denote this radius by $\rho(A)$
and set $t(A):=\cos\rho(A)$. The following is one of many equivalent
ways \cite{ChC00,ChCH05} of defining
the GCC condition number.
\begin{definition}
The GCC condition number of $A\in(S^m)^{n}$
is defined as ${\mathscr C}(A) := 1/|\cos\rho(A)| \in(1,\infty]$.
\end{definition}

In order to understand the relation of this condition number to
distance to ill-posedness, we review a few known facts (for more
information, see \cite{ChC00} and \cite{ChCH05}). Recall the definition
of $\mathcal{F}_{n,m}$ and $\mathcal {F}_{n,m}^\circ$ given in
equations (\ref{eq:F}) and (\ref{eq:Fcirc}). It is easy to see that
$\mathcal{F}_{n,m}$ is a compact semialgebraic set with nonempty
interior $\mathcal{F}_{n,m}^{\circ}$. The set
$\Sigma_{n,m}:=\mathcal{F}_{n,m}\setminus\mathcal {F}_{n,m}^\circ$ is
the topological boundary of $\mathcal{F}_{n,m}$. It consists of the
feasible instances that are not strictly feasible. Note that if
$n>m+1$, then $\Sigma_{n,m}$ is also the boundary of the set of
infeasible instances $\mathcal{I}_{n,m}:= (S^{m})^n
\setminus\mathcal{F}_{n,m}$. Hence in this case $\Sigma_{n,m}$ consists
of those instances that can be made both feasible and infeasible by
arbitrarily small perturbations.

The next lemma summarizes results from \cite{ChC00}, Theorem 1, and
\cite{ChCH05}, Proposition~4.1. It is enough to prove Proposition
\ref{prop:equiv}.
\begin{lemma}\label{le:hauser1}
We have $\rho(A)<\pi/2$ if and only if $A\in\mathcal
{F}_{n,m}^{\circ}$.
Moreover, $\rho(A)=\pi/2$ if and only if $A\in\Sigma_{n,m}$.
\end{lemma}
\begin{pf*}{Proof of Proposition \protect\ref{prop:equiv}}
We claim that
%
%
\begin{equation}\label{eq:cov-rho}
p(n,m,\alpha) = \mathsf{Prob}\{\rho(A)\le\pi-\alpha\}.
\end{equation}
Indeed, $\bigcup_{i=1}^n \mathsf{cap}(a_i,\alpha)\neq S^m$ iff there
exists $y\in S^m$ such that $y\notin\mathsf{cap}(a_i,\alpha)$ for all
$i$. This is equivalent to $\exists y\ \forall i\
a_i\notin\mathsf{cap}(y,\alpha)$ which means that an LEC for $A$ has
angular radius at least $\alpha$. This in turn is equivalent to
$\rho(A)\le\pi-\alpha$ thus proving the claim.

Equation (\ref{eq:cov-rho}) for $\alpha=\pi/2$ combined with Lemma
\ref{le:hauser1} and Wendel's result \cite{Wendel62} stated in equation
(\ref{eq:wendel}) yields
\[
2^{1-n}\sum_{k=0}^{m}\pmatrix{n-1\cr k}
= p(n,m,\pi/2)=\mathsf{Prob}\{\rho(A)\le\pi/2\} =
\mathsf{Prob}\{A\in\mathcal{F}_{n,m}\}.
\]
Suppose now $\alpha\in[\pi/2,\pi]$. Then
\[
\rho(A)\le\pi-\alpha\quad\iff\quad
\rho(A)\le\pi/2 \quad\mbox{and}\quad{\mathscr C}(A)\leq\frac1\varepsilon,
\]
showing the first assertion of Proposition \ref{prop:equiv}.
Furthermore, for $\alpha\in[0,\frac{\pi}{2}]$
\[
\rho(A)\le\pi-\alpha
\]
iff
\[
\rho(A)\le\pi/2\quad\mbox{or}\quad
\bigl(\rho(A) > \pi/2 \mbox{ and }
|{\cos\rho(A)}|\leq|{\cos(\pi-\alpha)}|\bigr),
\]
showing the second assertion of Proposition \ref{prop:equiv}.
\end{pf*}

\subsection[Toward the proof of Theorem 1.3]{Toward the proof of Theorem \protect\ref{th:main_lp}}
\label{sec:main_ideas}

To prove the feasible case in Theorem~\ref{th:main_lp}
we note that
\[
\mathsf{Prob}\biggl\{{\mathscr C}(A)\geq\frac1\varepsilon\Bigm|
A\in\mathcal{F}_{n,m} \biggr\}
=\frac{1}{\operatorname{\mathsf{vol}}\mathcal{F}_{n,m}}\operatorname{\mathsf{vol}}\mathcal
{F}_{n,m}(\varepsilon),
\]
where $\mathcal{F}_{n,m}(\varepsilon)=\{A\in\mathcal
{F}_{n,m}^\circ\mid t(A)< \varepsilon\}$. But
$\operatorname{\mathsf{vol}}\mathcal{F}_{n,m}$ is known by
Proposition~\ref{prop:equiv}. Therefore, our task is reduced to computing
$\operatorname{\mathsf{vol}}\mathcal{F}_{n,m}(\varepsilon)$. As we will
see in Section~\ref{se:sicprop}, the smallest including cap
$\mathrm{SIC}(A)$ is uniquely determined for all
$A\in\mathcal{F}_{n,m}^\circ$. Furthermore, for such $A$, $t(A)$
depends only on the blocking set of~$A$. Restricting to a suitable open
dense subset $\mathcal{R}_{n,m}(\varepsilon)\subseteq\mathcal
{F}_{n,m}(\varepsilon)$ of ``regular'' instances, these blocking sets
are of cardinality at most $m+1$. This induces a partition
\[
\mathcal{R}_{n,m}(\varepsilon)=\bigcup_I \mathcal
{R}_{n,m}^I(\varepsilon),
\]
where the union is over all subsets $I\subseteq[n]$ of cardinality at
most $m+1$, and $\mathcal{R}_{n,m}^I(\varepsilon)$ denotes the set of
matrices in $\mathcal{R}_{n,m}(\varepsilon)$ with blocking set indexed
by $I$. By symmetry,
$\operatorname{\mathsf{vol}}\mathcal{R}_{n,m}^I(\varepsilon)$ only
depends on the cardinality of $I$; hence it is enough to focus on
computing
$\operatorname{\mathsf{vol}}\mathcal{R}_{n,m}^{[k+1]}(\varepsilon)$ for
$k=1,\ldots,m$. The orthogonal invariance and the particular structure
of the $\mathcal{R}_{n,m}^{[k+1]}(\varepsilon)$ (involving certain
convexity conditions) makes possible a change of coordinates that
allows one to split the occurring integral into an integral over $t$
and a quantity $C(m,k)$ that depends only on $m$ and~$k$:
\[
\operatorname{\mathsf{vol}}\mathcal{R}_{n,m}^{[k+1]}(\varepsilon)
=C(m,k)\int_{0}^{\varepsilon} g(t,n,m,k)\,dt.
\]

More precisely, we proceed as follows:

(1)
By Fubini, we split the integral over $\mathcal
{R}_{n,m}^{[k+1]}(\varepsilon)$
into an integral over the first $k+1$ vectors
$a_1,\ldots,a_{k+1}$
(determining the blocking set $[k+1]$)
and an integral over
$a_{k+2},\ldots,a_n$ taken from $\mathrm{SIC}(A)$:
%
%
\begin{eqnarray} \label{eq:volpow}\qquad
\operatorname{\mathsf{vol}}\mathcal{R}_{n,m}^{[k+1]}(\varepsilon)
&=&\int_{A\in\mathcal{R}_{k+1,m}^{[k+1]}(\varepsilon)}
\biggl(\int_{\mathsf{cap}(p(A),\rho(A))^{n-k-1}}d(S^m)^{n-k-1} \biggr)
\,d\mathcal{R}_{k+1,m}^{[k+1]} \nonumber\\[-8pt]\\[-8pt]
&=&\int_{A\in\mathcal{R}_{k+1,m}^{[k+1]}(\varepsilon)} G(A)\,
d\mathcal{R}_{k+1,m}^{[k+1]}.\nonumber
\end{eqnarray}
This is an integral of the function
$G(A) := \operatorname{\mathsf{vol}}(\mathsf{cap}(p(A),\rho(A)))^{n-k-1}$ which
is a certain power of the volume of the spherical cap $\mathrm{SIC}(A)$.

\begin{figure}

\includegraphics{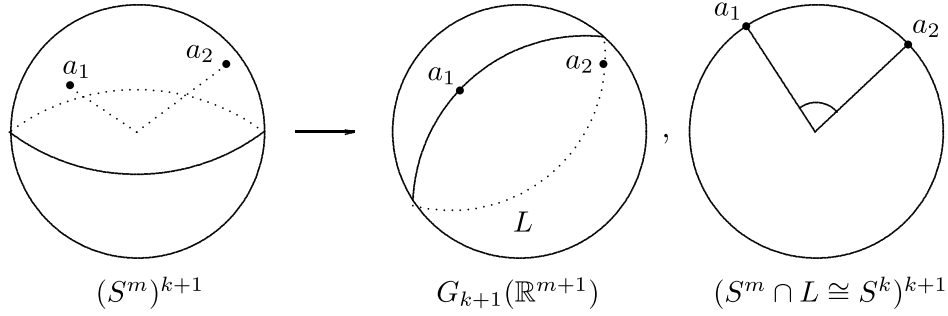}

\caption{Determining $(a_1,a_2)\in(S^2)^{2}$ by first giving its
span $L\in\mathrm{G}_{2}(\mathbb R^{3})$ and then $a_i\in S^2\cap L$.}
\label{fig:idea1}
\end{figure}

(2)
The next idea is to specify $A=(a_1,\ldots,a_{k+1})$ in $\mathcal{R}
_{k+1,m}^{[k+1]}(\varepsilon)$
by first specifying the subspace $L$ spanned by these vectors
and then the position of the $a_i$
on the sphere $S^m\cap L\cong S^{k}$ (cf. Figure \ref{fig:idea1}).

%
%

Let $\mathrm{G}_{k+1}(\mathbb R^{m+1})$ denote the Grassmann manifold of
$(k+1)$-dimensional subspaces in $\mathbb R^{m+1}$ and
consider the map
%
%
\begin{eqnarray}\label{eq:mapuno}
\mathcal{R}_{k+1,m}^{[k+1]}(\varepsilon)&\rightarrow&\mathrm
{G}_{k+1}(\mathbb R^{m+1}),\nonumber\\[-8pt]\\[-8pt]
(a_1,\ldots,a_{k+1})&\mapsto&L=\mathsf{span}\{a_1,\ldots,a_{k+1}\}.
\nonumber
\end{eqnarray}
Clearly, a vector $a\in S^m$ lies in the special subspace $L_0:=\mathbb
R^{k+1}\times0$ iff it lies in the subsphere $S^k$. Hence the fibre
over $L_0$ consists of all ``regular'' tuples $A=(a_1,\ldots,a_{k+1})
\in(S^k)^{k+1}$ such that $t(A)<\varepsilon$, and hence the fibre can
be identified with $\mathcal{R}_{k+1,k}^{[k+1]}(\varepsilon)$. Using
the orthogonal invariance and the coarea formula (also called Fubini's
theorem for Riemannian manifolds) we can reduce the computation of the
integral (\ref{eq:volpow}) to an integral over the special fibre
$\mathcal {R}_{k+1,k}^{[k+1]}(\varepsilon)$. This leads to
\[
\int_{\mathcal{R}_{k+1,m}^{[k+1]}(\varepsilon)}
G(A) \,d\mathcal{R}_{k+1,m}^{[k+1]}=
\operatorname{\mathsf{vol}}\mathrm{G}_{k+1}(\mathbb R^{m+1})
\int_{\mathcal{R}_{k+1,k}^{[k+1]}(\varepsilon)} G(A)J(A) \,d\mathcal
{R}_{k+1,k}^{[k+1]},
\]
where $J(A)$ is the normal Jacobian of the transformation (\ref{eq:mapuno}).

\begin{figure}

\includegraphics{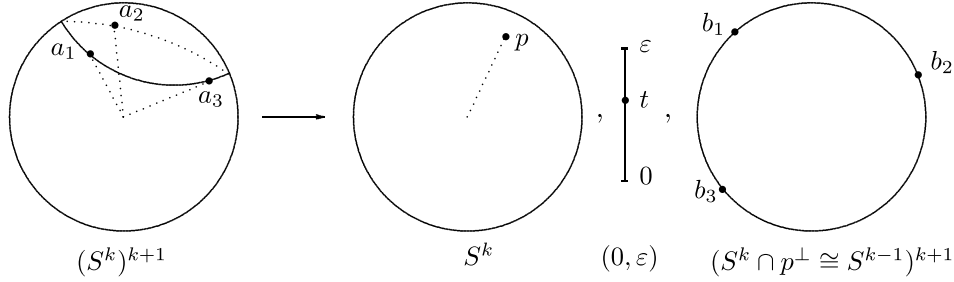}

\caption{Determining $(a_1,a_2,a_3)\in(S^2)^{3}$ by specifying
the direction $p$, the height $t$, and the $a_i$ on the subsphere
$\{a\in S^2\mid\langle a,p\rangle=t\}$ by $b_i$.}
\label{fig:idea2}
\end{figure}

(3)
To specify a regular $A=(a_1,\ldots,a_{k+1})\in(S^k)^{k+1}$,
we first specify the direction
$p=p(A)\in S^k$ and the height $t=t(A)\in(0,\varepsilon)$
and then the position of the $a_i$ on the subsphere
$\{a\in S^k\mid\langle a,p\rangle=t\}\simeq S^{k-1}$
(cf. Figure \ref{fig:idea2}).

%
%

More precisely, we consider the map
%
%
\begin{equation}\label{eq:mapdos}
\mathcal{R}_{k+1,k}^{[k+1]}(\varepsilon)\rightarrow S^k\times
(0,\varepsilon),\qquad
A \mapsto(p(A),t(A)).
\end{equation}
The fibre over $(p_0,t)$, where $p_0=(0,\ldots,0,1)$ is the ``north
pole,'' consists of tuples $(a_1,\ldots,a_{k+1})$ lying on the
``parallel'' subsphere $\{a\in S^k\mid\langle a,p\rangle=t\}$. The
vectors $a_i$ can be described by points $b_i\in S^{k-1}$, which are
obtained by projecting $a_i$ orthogonally onto $\mathbb R^k\times0$ and
scaling the resulting vector to length one.

The orthogonal invariance and the coarea formula allow us to reduce the
computation of the integral over $\mathcal
{R}_{k+1,k}^{[k+1]}(\varepsilon)$ to the integration over
$t\in[0,\varepsilon]$ of an integral over the special fibres over
$(p_0,t)$. The latter integral is captured by the coefficient $C(m,k)$
which can be interpreted as a higher moment of the volume of the
simplex $\Delta$ spanned by random points $b_1,\ldots,b_{k+1}$ on the
sphere $S^{k-1}$. However, we have to respect the convexity condition
that the origin is contained in the simplex $\Delta$ spanned by the
$b_i$, which complicates matters. Altogether, we are lead to a formula
for $\operatorname{\mathsf{vol}}\mathcal{R}_{n,m}^{[k+1]}(\varepsilon)$
of the shape,
\begin{eqnarray*}
&&\operatorname{\mathsf{vol}}\mathrm{G}_{k+1}(\mathbb R^{m+1})
\int_{\mathcal{R}_{k+1,k}^{[k+1]}(\varepsilon)} G(A)J(A)\, d\mathcal
{R}_{k+1,k}^{[k+1]}
\\
&& \qquad= C(m,k) \int_{0}^{\varepsilon} g(t,n,m,k) \, dt,
\end{eqnarray*}
where $g(t,n,m,k)$ is obtained by isolating the part of the resulting
integrand that depends on $t$.

In order to implement this plan, we have to
isolate the appropriate regularity conditions, that is,
to identify the sets $\mathcal{R}_{n,m}^{I}(\varepsilon)$,
and to compute the normal Jacobians of the
maps (\ref{eq:mapuno}) and (\ref{eq:mapdos}).
For the latter task, we prefer to use the language of differential forms
as is common in integral geometry \cite{sant76}.

Unfortunately, the above argument does not carry over to the
infeasible case. Nevertheless, the ideas described above
are sufficient to obtain the upper bound in Theorem \ref{th:main_lp}.

The rest of the paper proceeds as follows. In Section \ref{se:prelim}
we describe the basic facts on smallest including caps and integration
on manifolds that will be needed to make formal the ideas expressed
above. Then, in Section \ref{sec:1.3}, we prove
Theorem~\ref{th:main_lp}. Theorem \ref{th:main_coverage} immediately
follows via Proposition \ref{prop:equiv}. Finally, in
Section~\ref{se:C(m,k)}, we give bounds for all, explicit expressions for some,
and a way to compute the coefficients $C(m,k)$. From these bounds we
derive Theorems \ref{cor:2} and \ref{corol:1}.

\section{Preliminaries}\label{se:prelim}

\subsection{Properties of smallest including caps}\label{se:sicprop}

Recall from Section \ref{se:prel-poly} the definition of smallest
including caps (SICs) for a given $A=(a_1,\ldots,a_n)\in(S^m)^n$.
A~crucial feature of our proofs is the fact that a strictly feasible
$A$ has a uniquely determined SIC. This is a consequence of the
following crucial lemma which provides an explicit criterion for a
spherical cap being a smallest including cap of $A$. This lemma is a
generalization of Lemma 4.5 in \cite{ChCH05}.
\begin{lemma}\label{le:crit}
\textup{(a)} For a strictly feasible $A\in\mathcal{F}_{n,m}^{\circ}$ there
exists exactly one smallest including cap.\smallskipamount=0pt

\begin{longlist}[(b)]
\item[(b)] Let $(p,t)\in S^m\times(0,1]$ and $1\le k < n$.
Suppose that
$\langle a_i,p\rangle= t$ for all $i\in[k+1]$ and
$\langle a_i,p\rangle> t$ for all $ i\in[n]\setminus[k+1]$.
Then $\mathsf{cap}(p,\arccos t)$ is
the smallest including cap of $A$ if and only if
\[
tp \in\mathsf{conv}\{a_1,\ldots,a_{k+1}\} .
\]
\end{longlist}
\end{lemma}
\begin{pf}
We first show that assertion (b) implies assertion (a).
Suppose $\mathsf{cap}(p_1,\rho)$ and $\mathsf{cap}(p_2,\rho)$ are
SICs for $A$,
and put $t:= \cos\rho$. Note that $t>0$.
Assertion (b) implies that $tp_1$ is contained
in the convex hull of $a_1,\ldots,a_n$; hence there exist $\lambda
_i\ge0$
such that $\sum_i\lambda_i=1$ and
$tp_1=\sum_i\lambda_i a_i$. Therefore,
$\langle tp_1,p_2 \rangle= \sum_i \lambda_i \langle a_i,p_2 \rangle
\ge t$,
as $\langle a_i,p_2 \rangle\ge t$ for all $i$.
This implies $\langle p_1,p_2\rangle\ge1$ and hence $p_1=p_2$.

The proof of assertion (b) goes along the lines of Lemma 4.5 in \cite{ChCH05}.
Suppose first that $\mathsf{cap}(p,\alpha)$ is a SIC for $A$ where
$\alpha
:=\arccos t$.
It is sufficient to show that $p\in\mathsf{cone}\{a_1,\ldots
,a_{k+1}\}$.
Indeed, if $p=\sum\lambda_ia_i$ with $\lambda_i\ge0$, then
$tp=\sum(t\lambda_i)a_i$. Furthermore,
$\sum(t\lambda_i) =\sum\lambda_i\langle a_i,p\rangle
= \langle\sum\lambda_i a_i,p \rangle=\|p\|^2=1$.
Hence
$tp\in\mathsf{conv}\{a_1,\ldots,a_{k+1}\}$.

We now argue by contradiction:
if $p$ is not contained in $\mathsf{cone}\{a_1,\ldots,a_{k+1}\}$, then
there exists a vector $v\in S^m$ such that $\langle p,v\rangle<0$ and
$\langle a_i,v\rangle>0$ for all $i\in\{1,\ldots,k+1\}$.
For $\delta>0$ we set
%
%
\begin{equation}\label{eq:pdelta}
p_{\delta}:=\frac{p+\delta v}{\|p+\delta v\|}
= \frac{p+\delta v}{\sqrt{1+2\delta\langle p,v\rangle+\delta^2}}.
\end{equation}
Then for $1\leq i\leq k+1$ and sufficiently small $\delta$ we have
\[
\langle a_i,p_{\delta}\rangle
=\frac{t+\delta\langle a_i,v\rangle}{\sqrt{1+2\delta\langle
p,v\rangle+\delta^2}}>t,
\]
where we used that $\langle a_i,p\rangle=t$, $\langle a_i,v\rangle
>0$ and $\langle p,v\rangle<0$.

For $k+2\leq i\leq n$ the function $\delta\rightarrow\langle
a_i,p_{\delta}\rangle$
is continuous at $\delta=0$. Since, by hypothesis,
$\langle a_i,p\rangle=\langle a_i,p_0\rangle>t$, it follows that
$\langle a_i,p_\delta\rangle> t$
for $\delta$ sufficiently small. From this we conclude that for
sufficiently small $\delta$ there
exists $t_\delta>t$ such that $\langle a_i,p_\delta\rangle>t_\delta
$ for all
$i\in[n]$. Setting $\alpha_\delta=\arccos t_\delta$ we obtain that
$\alpha_\delta<\alpha$ and $a_i\in
\mathsf{cap}(p_\delta,\alpha_\delta)$ for all $i\in[n]$, contradicting
the assumption that $\mathsf{cap}(p,\alpha)$ is a smallest including cap.

To prove the other direction, we suppose $tp\in\mathsf{conv}\{
a_1,\ldots ,a_{k+1}\}$. For $q\in S^m$ let $\alpha(q)$ we denote the
angular radius of the smallest spherical cap with center $q$ containing
$a_1,\ldots, a_n$. If we assume that $\mathsf{cap}(p,\alpha)$ is not a
SIC for $A$, then there exists a vector $v\in S^m$ and $\delta_0>0$,
such that $\langle v,p\rangle =0$ and, for all $0<\delta\leq\delta_0$,
$\alpha(p_\delta)<\alpha(p)$ where $p_\delta=\frac{p+\delta
v}{\sqrt{1+\delta^2}}$ (i.e., we have a direction $v$ along which we
can move to obtain a smaller cap). This means that
\[
\min_{1\leq i\leq n} \langle a_i,p_\delta\rangle>\min_{1\leq i\leq
n}\langle a_i,p\rangle=t.
\]
Therefore, for all $i\in[k+1]$ we have
\[
\langle a_i,p_\delta\rangle=\frac{\langle a_i,p\rangle+\delta
\langle a_i,v\rangle}{\sqrt{1+\delta^2}}
> t=\langle a_i,p\rangle
\]
for sufficiently small $\delta$
which implies that $\langle a_i,v\rangle> 0$.
Let $\mu\in\mathbb R^{k+1}_{\geq0}$ be such that
$tp=\sum_{1\leq i\leq k+1}\mu_ia_i$ and $\sum_{1\leq i\leq k+1}\mu_i=1$.
Then we have
\[
t\langle p,v\rangle=\sum_{1\leq i\leq k+1}\mu_i\langle a_i,v\rangle>0,
\]
contradicting the assumption that $\langle p,v\rangle=0$.
Thus $\mathsf{cap}(p,\alpha)$ is indeed a smallest including cap.
\end{pf}

For a strictly feasible $A$, we denote the center of its uniquely
determined SIC by $p(A)$ and its radius by $\rho(A)$. The blocking set
$\mathrm{BS}(A)$ of $A$ is defined as the blocking set of the SIC of
$A$. It is not hard to see that $\mathrm{BS}(A)$ can have any
cardinality greater than one.

However, we note that in the infeasible case, there may be more than
one smallest including cap. Consider for instance three equilateral
points on the circle (right-hand side in Figure \ref{fig:1}). It is
known \cite{ChCH05}, Proposition 4.2, that in this case, the blocking
set of a SIC has at least $m+1$ elements. In the infeasible case, one
direction of the characterization of smallest including caps of Lemma
\ref{le:crit} still holds. The proof is similar as for Lemma
\ref{le:crit}.
\begin{lemma}\label{lem:new}
Let $\mathsf{cap}(p,\arccos t)$ be a $\mathrm{SIC}$ for $A\in(S^m)^n$
with $p\in S^m$ and $t\in(-1,0)$. Suppose $\langle a_i,p\rangle=t$ for
$i\in[m+1]$ and $\langle a_i,p\rangle>t$ for $i=m+2,\ldots,n$. Then
$tp\in\mathsf{conv}\{a_1,\ldots,a_{m+1}\}$.
\end{lemma}
\begin{pf}
Suppose $tp\notin\mathsf{conv}\{a_1,\ldots,a_{m+1}\}$. Then
$-p\notin\mathsf{cone}\{a_1,\ldots,a_{m+1}\}$ and hence there exists a
vector $v\in S^m$ such that $\langle-p,v\rangle <0$ and $\langle
a_i,v\rangle>0$ for all $i$. For $\delta>0$ we define $p_\delta$ as in
(\ref{eq:pdelta}). Take $\delta$ sufficiently small so that
$t<t+\delta\langle a_i,v\rangle<0$ for all $i\in[m+1]$. Then, for
$i\in[m+1]$ and $\delta$ sufficiently small, we have
\[
\langle a_i,p_{\delta}\rangle =\frac{t+\delta\langle
a_i,v\rangle}{\sqrt{1+2\delta\langle p,v\rangle+\delta^2}}>t,
\]
where we used that $\langle a_i,p\rangle=t$, $\langle a_i,v\rangle
>0$, and $\langle p,v\rangle>0$.
This shows that $\mathsf{cap}(p,\arccos t)$ is not a smallest
including cap.
\end{pf}

We present a few more auxiliary results that are needed for the proof
of our main result.
\begin{lemma}\label{le:pt}
For given linearly independent $a_1,\ldots,a_{k+1}\in S^m$, $1\leq
k\leq m$, there exist \textit{unique} $p\in S^m$ and $t\in(0,1)$ such
that
\[
p\in\mathsf{span}\{a_1,\ldots,a_{k+1}\}
\]
and
\[
\langle a_i,p\rangle= t \qquad\mbox{for all } i\in[k+1].
\]
Moreover, the map $(a_1,\ldots,a_{k+1})\mapsto(p,t)$ is differentiable.
\end{lemma}
\begin{pf}
Let $\mathcal{A}$ denote the affine span of $a_1,\ldots,a_{k+1}$, $L$
the underlying linear space and $\mathcal{L}$ the linear span of
$\mathcal{A}$. Since the $a_i$ are linearly independent, we have
$\mathcal{A}\ne L$ and thus $\dim\mathcal{A}=\dim L=k$,
$\dim\mathcal{L}=k+1$. Hence the intersection of $\mathcal{L}$ with the
orthogonal complement $L^{\perp}$ is one-dimensional and contains
exactly two elements of length one. Take the one such that the common
value $t=\langle a_i,p\rangle$ is positive. This shows existence and at
the same time the uniqueness of $p,t$.

Suppose now $k=m$, and let $A$ denote the square matrix with the rows
$a_1,\ldots,a_{m+1}$. The conditions $\langle a_i,p\rangle= t$ can be
written in matrix form as $Ap=te$ where $e:=(1,\ldots,1)^{\top}$. By
solving this equation we obtain the following explicit formulas:
%
%
\begin{equation}\label{eq:pt-formula}
p(A) = \frac1{\|A^{-1}e\|} A^{-1}e,\qquad t(A) = \frac1{\|A^{-1}e\|} .
\end{equation}
This shows the differentiability of the map $A\mapsto(p,t)$ in the case
$k=m$. We leave the proof in the general case to the reader.
\end{pf}

The next result, though very elementary, will be useful for
clarification.

Let $p\in S^{k}$ and $t\ne0$ and consider elements
$a_1,\ldots,a_{k+1}\in S^{k}$ satisfying $\langle a_i,p\rangle= t$ for
all $i$. Let $b_i\in S^{k-1}$ be the scaled-to-one orthogonal
projection of $a_i$ onto the orthogonal complement of $\mathbb Rp$.
That is, $a_i=rb_i +tp$ where $r=(1-t^2)^{1/2}$.
\begin{lemma}\label{le:elusive}
The following conditions are equivalent:
\begin{enumerate}
\item The affine hull $\mathcal{A}$ of $a_1,\ldots,a_{k+1}$ has
dimension $k$.
\item The span of $b_1,\ldots,b_{k+1}$ has dimension $k$.
\item$a_1,\ldots,a_{k+1}$ are linearly independent.
\end{enumerate}
\end{lemma}
\begin{pf}
The equivalence of the first two conditions is obvious.
The equivalence of the first and third condition follows from
$\dim\mathsf{span}(\mathcal{A})= \dim\mathcal{A} +1$
(here we use $t\ne0$).
\end{pf}

\subsection{Volume forms on Grassmann manifolds}\label{se:Grass}

Integration on Grassmann mani\-folds will play a crucial role in our proofs.
We recall some facts about the relevant techniques from integral geometry
and refer to Santal\'{o}'s book \cite{sant76}, II.9--12, and
the article \cite{mile71} for more information.
We recall that volume elements are always unsigned forms.

Let $M$ be a Riemannian manifold of dimension $m$, $p\in M$, and let
$y=(y_1,\ldots,y_m)^{\top}\dvtx U\mapsto\mathbb R^m$ be local
coordinates in a neighborhood $U$ of $p$ such that $\partial/\partial
y_1,\ldots ,\partial/\partial y_m$ are an orthonormal basis of $T_pM$.
The natural volume form on $M$ at $p$ associated to its Riemannian
metric is then given by $dM= dy_1\wedge\cdots\wedge dy_m$ where $dy_i$
is the differential of the coordinate function $y_i$ at $p$.

In the case of a sphere, we get such local coordinates around a point
$p\in S^m$ by projecting onto the orthogonal complement of $p$. More
precisely, let $\langle e_1,\ldots,e_{m+1}\rangle$ be an orthonormal
basis of $\mathbb R^{m+1}$ satisfying $e_1=p$ (so that
$e_2,\ldots,e_{m+1}$ span the tangent space $T_pS^m$). For a point
$x=(x_1,\ldots,x_{m+1})^{\top}\in S^m$ in a neighbourhood of $p$ set
$y_i=\langle x,e_i\rangle$. Then $(y_2,\ldots,\break y_{m+1})$ are local
coordinates of $S^m$ around $p$ such that $\partial/\partial y_i$ are
pairwise orthogonal at $p$. Hence the volume element of $S^m$ at $p$ is
given by
\[
dS^m = \omega_2\wedge\cdots\wedge\omega_{m+1},
\]
where $\omega_i:=dy_i=\langle dx,e_i\rangle$ and $dx=(dx_1,\ldots
,dx_{m+1})^{\top}$.
Hence, if we denote by $E$ the $(m+1)\times(m+1)$-matrix
having the $e_i$ as rows,
we obtain the volume form by wedging the
nonzero entries of $E\, dx$.

In a similar fashion we define volume forms on Stiefel manifolds (for
details and further justification we refer to \cite{sant76}). A
$k$-\textit{frame} is a set of $k$ linearly independent vectors. For
$1\leq k\leq m+1$, the \textit{Stiefel manifold} $V_k(\mathbb R^{m+1})$
is defined as the set of orthonormal $k$-frames in $\mathbb R^{m+1}$.
It is a compact Riemannian submanifold of $(S^{m})^k$. Let
$Q=(q_1,\ldots,q_k)\in V_{k}(\mathbb R^{m+1})$ and $\langle e_1,\ldots
,e_{m+1}\rangle$ an orthonormal basis of $\mathbb R^{m+1}$ such that
$e_1=q_1,\ldots,e_k=q_k$. Then the volume element of $V_{k}(\mathbb
R^{m+1})$ at $Q$ is given by
\[
dV_{k}(\mathbb R^{m+1})=\bigwedge_{1\leq i\leq k}
( \omega_{i,i+1}\wedge\cdots\wedge\omega_{i, m+1} ),
\]
where $\omega_{i,j}=\langle dq_i,e_j\rangle$ for $1\leq i\leq k$ and
$1\leq j\leq m+1$. [In terms of the $(m+1)\times k$ matrix $E \,dQ$, this
corresponds to wedging the entries below the main diagonal.] With this
volume element we have $\operatorname{\mathsf{vol}}V_k(\mathbb
R^{m+1})=\mathcal O_m \cdots\mathcal O_{m+1-k}$.

We denote by $\mathrm{G}_k(\mathbb R^{m+1})$ the \textit{Grassmann
manifold} of $k$-dimensional subspaces of $\mathbb R^{m+1}$. One way of
characterizing it is as a quotient of a Stiefel manifold, by
identifying frames that span the same subspace. Let
$L\in\mathrm{G}_k(\mathbb R^{m+1})$ and choose a frame $Q\in
V_{k}(\mathbb R^{m+1})$ spanning $L$. If $V_k(L)$ denotes the Stiefel
manifold of orthonormal $k$-frames in $L$ and $dV_k(L)$ its volume
element at $Q$, then it is known that the volume element
$d\mathrm{G}_{k}(\mathbb R^{m+1})$ of the Grassmann manifold at $L$
satisfies (see \cite{mile71}, equation (10))
%
%
\begin{equation}\label{eq:grass-stiefel}
dV_{k}(\mathbb R^{m+1}) = d\mathrm{G}_{k}(\mathbb R^{m+1})\wedge dV_k(L).
\end{equation}
From this equality it follows that
\[
d\mathrm{G}_k(\mathbb R^{m+1})=\bigwedge_{1\leq i\leq k}
(\omega_{i,k+1}\wedge\cdots\wedge\omega_{i, m+1} )
\]
with the $\omega_{i,j}$ as defined in the case of the Stiefel manifold.
(In terms of the matrix $E dQ$, this corresponds to wedging the
elements in the lower $(m+1-k)\times k$ rectangle.) As a consequence of
(\ref{eq:grass-stiefel}), the volume of the Grassmannian is given by
\[
\mathcal G_{k,m+1}:=\operatorname{\mathsf{vol}}\mathrm{G}_k(\mathbb R^{m+1})
=\frac{\mathcal O_{m+1-k}\cdots\mathcal O_{m}}{\mathcal O_0\cdots
\mathcal O_{k-1}}.
\]

Equation (\ref{eq:grass-stiefel}) has a generalization to frames that
are not orthogonal, that is, to points in a product of spheres
$(S^{m})^{k}$. Let $L\in\mathrm{G}_{k}(\mathbb R^{m+1})$ and set
$S(L):=L\cap S^m$, so that $S(L)\cong S^{k-1}$. Choose a basis
$a_1,\ldots,a_{k}$ of $L$ consisting of unit length vectors, that is, a
point $A=(a_1,\ldots,a_{k})$ in $S(L)^{k}$. We denote by
$\operatorname{\mathsf{vol}}(A)$ the volume of the parallelepiped
spanned by the $a_i$. Then the volume form of $(S^m)^{k}$ at $A$ can be
expressed as
%
%
\begin{equation}\label{eq:useful-trans}
d(S^m)^{k}= \operatorname{\mathsf{vol}}(A)^{m-k+1} \, d\mathrm{G}_{k}(\mathbb
R^{m+1})\wedge dS(L)^{k}.
\end{equation}
This equation can be derived as in \cite{mile71} (see also \cite
{sant76}, II.12.3). It also follows as a special case of a general
formula of Blaschke--Petkantschin-type derived by Z\"{a}hle
\cite{zaehle} (see also the discussion in \cite{reitz}).

A beautiful application of equation (\ref{eq:useful-trans}) is that it
allows an easy computation of the moments of the absolute values of
random determinants. The following lemma is an immediate consequence
of (\ref{eq:useful-trans}) (see also \cite{mile71}).
\begin{lemma}\label{le:highmom}
Let $B\in(S^{k})^{k+1}$ be a matrix with rows $b_1,\ldots,b_{k}$
independently and uniformly distributed in $S^{k}$. Then
\[
\mathbf{E}(|{\det(B)}|^{m-k+1})
= \biggl(\frac{\mathcal O_{m}}{\mathcal O_{k-1}} \biggr)^{k}
\frac{1}{\mathcal G_{k,m+1}}.
\]
\end{lemma}

\section{The probability distribution of ${\mathscr{C}}(A)$}\label{sec:1.3}

This section is devoted to the proof of Theorem \ref{th:main_lp}.

\subsection{The feasible case}\label{se:nJac}

Recall that, for $A\in\mathcal{F}_{n,m}^{\circ}$, we
denote the center and the angular radius of the unique smallest
including cap of $A$ by $p(A)$ and $\rho(A)$, respectively,
and we write $t(A)=\cos\rho(A)$.

Our goal here is to prove the first part of Theorem \ref{th:main_lp},
for which, as we noted in Section \ref{sec:main_ideas}, we just need to
compute the volume of the following sets:
\[
\mathcal{F}_{n,m}(\varepsilon): =\{A\in\mathcal{F}_{n,m}^\circ
\mid t(A)<\varepsilon\}.
\]
For this purpose
it will be convenient to decompose
$\mathcal{F}_{n,m}(\varepsilon)$ according to the size of the
blocking sets.
Recall that the \textit{blocking set} of $A\in\mathcal{F}_{n,m}^\circ$
is defined as
%
%
\begin{equation}\label{eq:defBS}
\mathrm{BS}(A)=\{i\in[n]\mid\langle p(A),a_i\rangle= t(A) \} .
\end{equation}
For $I\subseteq[n]$ with $|I|\le n$ and $\varepsilon\in(0,1]$ we define
$\mathcal{F}_{n,m}^I(\varepsilon)$ to be the set of all $A\in
\mathcal{F}_{n,m}(\varepsilon)$
such that $\mathrm{BS}(A)=I$.

For technical reasons we have to require some regularity conditions
for the elements of $\mathcal{F}_{n,m}^I(\varepsilon)$.
\begin{definition}\label{def:centered}
We call a family $(a_1,\ldots,a_{k+1})$ of elements of a vector space
\textit{centered with respect to a vector $c$} in the affine hull
$\mathcal{A}$
of the $a_i$ if $\dim\mathcal{A} =k$, and
$c$ lies in the relative interior of the convex hull of the $a_i$.
We call the family \textit{centered} if it is centered with respect
to the origin. We now define, for $I\subseteq[n]$,
\[
\mathcal{R}^I_{n,m}(\varepsilon) := \{ A\in\mathcal
{F}_{n,m}^I(\varepsilon) \mid(a_i)_{i\in I}
\mbox{ is centered with respect to $t(A)p(A)$} \}.
\]
\end{definition}

Note that, by definition, the $a_i$ are affinely
independent if $|I|\leq m+1$.
\begin{lemma}\label{le:Fopen}
\textup{1. }
$\mathcal{F}_{n,m}^I(\varepsilon)$ is of measure zero if $|I|> m+1$.
\smallskipamount=0pt
\begin{enumerate}[2.]
\item[2.] If $|I|\le m+1$, then $\mathcal{R}_{n,m}^I(\varepsilon)$ is
open in $(S^m)^n$,
and $\mathcal{F}_{n,m}^I(\varepsilon)$ is contained in the closure
of $\mathcal{R}_{n,m}^I(\varepsilon)$.
\item[3.] $\mathcal{F}_{n,m}^I(\varepsilon)\setminus\mathcal
{R}_{n,m}^I(\varepsilon)$ has measure zero.
\end{enumerate}
\end{lemma}
\begin{pf}
1. If $A\in\mathcal{F}_{n,m}^I(\varepsilon)$, then $\{a_i\mid i\in
I\}$ is contained in
the boundary of the SIC of $A$, and hence its affine hull has dimension
at most $m$. On the other hand, the affine hull of $(a_i)_{i\in I}$
is almost surely $\mathbb R^{m+1}$ if $|I|> m+1$.

2. The fact that $\mathcal{R}_{n,m}^I(\varepsilon)$ is open in
$(S^m)^n$ easily follows
from the continuity of the map $A\mapsto(p(A),t(A))$ established in
Lemma \ref{le:pt}.

Suppose now $A\in\mathcal{F}_{n,m}^I(\varepsilon)$.
By Lemma \ref{le:crit} we have
$t(A)p(A)\in\mathsf{conv}\{a_i \mid i\in I\}$ for $A\in\mathcal
{F}_{n,m}^I(\varepsilon)$.
It is now easy to see that there are elements $A'$
arbitrarily close to $A$ such that
$A'$ is centered with respect to $t(A')p(A')$.
This shows the second assertion.

3. By part two we have
$\mathcal{R}_{n,m}^I(\varepsilon)\subseteq\mathcal
{F}_{n,m}^I(\varepsilon)\subseteq\overline{\mathcal{R}
_{n,m}^I(\varepsilon)}$.
Since we are dealing with semialgebraic sets,
the boundary of $\mathcal{R}_{n,m}^I(\varepsilon)$ is of measure zero.
\end{pf}

It is clear that the $\mathcal{F}_{n,m}^I$
with $I$ of the same cardinality just differ
by a permutation of the vectors.
Using Lemma \ref{le:Fopen} we obtain
%
%
\begin{equation}\label{eq:partition}
\operatorname{\mathsf{vol}}\mathcal{F}_{n,m}(\varepsilon) = \sum_{|I|\le m+1}
\operatorname{\mathsf{vol}}\mathcal{F}_{n,m}^{I}(\varepsilon)
= \sum_{k=1}^m \pmatrix{n\cr{k+1}} \operatorname{\mathsf{vol}}\mathcal
{R}_{n,m}^{k}(\varepsilon),
\end{equation}
where we have put
$\mathcal{R}_{n,m}^{k}(\varepsilon) := \mathcal
{R}_{n,m}^{[k+1]}(\varepsilon)$ to ease notation.

Hence it is sufficient to compute the volume of $\mathcal
{R}_{n,m}^{k}(\varepsilon)$.
For this purpose we introduce now the coefficients $C(m,k)$.
\begin{definition}\label{def:C}
We define for $1\le k\le m$
\[
C(m,k)=\frac{(k!)^{m-k+1}}{\mathcal O_m^k} \mathcal G_{k,m}\int_{M_k}
(\operatorname{\mathsf{vol}}_{k} \Delta)^{m-k+1}\, d(S^{k-1})^{k+1},
\]
where $M_k$ is the following open subset of $S^{k-1}$:
\[
M_k := \{(b_1,\ldots,b_{k+1})\in(S^{k-1})^{k+1} \mid
\mbox{$(b_1,\ldots,b_{k+1})$ is centered}\}
\]
and $\Delta\dvtx M_k\to\mathbb R$ maps $B=(b_1,\ldots,b_{k+1})$ to the
convex hull of the $b_i$.
\end{definition}
\begin{example}\label{ex:k=2}
We compute $C(m,1)$.
Note that $M_{1,1}=\{(-1,1),(1,-1)\}$
and $\mathcal G_{1,m}=\frac12\mathcal O_{m-1}$. Hence
\[
C(m,1) = \frac{1}{\mathcal O_{m}} \frac12\mathcal O_{m-1}
\int_{M_{1,1}} (\operatorname{\mathsf{vol}}_{1}\Delta)^{m}\, dM_{1,1}
= \frac{\mathcal O_{m-1}}{\mathcal O_{m}} 2^{m}.
\]
\end{example}

The assertion in Theorem \ref{th:main_lp}
about the feasible case follows immediately from
Proposition \ref{prop:equiv}, equation (\ref{eq:partition})
and the following result.
\begin{proposition}\label{th:main-f}
Let $\varepsilon\in(0,1]$. The relative volume of $\mathcal
{F}_{n,m}^{k}(\varepsilon)$ is
given by
\[
\frac{\operatorname{\mathsf{vol}}\mathcal{R}_{n,m}^{k}(\varepsilon)}{\mathcal O_m^n}
= C(m,k)
\int_0^{\varepsilon} t^{m-k}(1-t^2)^{km/2-1}
\lambda_m(t)^{n-k-1} \,dt.
\]
\end{proposition}
\begin{pf}
Consider the projection
\[
\mathcal{R}_{n,m}^{k}(\varepsilon)\to\mathcal
{R}_{k+1,m}^{k}(\varepsilon),\qquad
(a_1,\ldots,a_n)\mapsto A=(a_1,\ldots,a_{k+1}) .
\]
By Lemma \ref{le:crit}, this map is surjective and
its fiber over $A$ consists of all
$(A,a_{k+2},\ldots,a_n)$ such that $a_i$ lies in the
interior of the cap $\mathsf{cap}(p(A),\rho(A))$ for all $i>k+1$.
By Fubini, and using (\ref{eq:volcap}), we conclude that
%
%
\begin{equation}\label{eq:weight}
\frac{\operatorname{\mathsf{vol}}\mathcal{R}_{n,m}^{k}(\varepsilon)}{\mathcal
O_m^{n-k-1}} = \int_{A\in\mathcal{R}_{k+1,m}^{k}(\varepsilon)} \lambda
_m(t(A))^{n-k-1} \,d(S^m)^{k+1}.
\end{equation}
We consider now the following map (which is well defined [cf. Lemma
\ref{le:elusive}]):
\[
\mathcal{R}_{k+1,m}^k(\varepsilon) \to\mathrm{G}_{k+1}(\mathbb
R^{m+1}),\qquad
(a_1,\ldots,a_{k+1})\mapsto L=\mathsf{span}\{a_1,\ldots,a_{k+1}\}.
\]
We can thus integrate over $\mathcal{R}_{k+1,m}^k(\varepsilon)$ by
first integrating over
$L\in\mathrm{G}_{k+1}(\mathbb R^{m+1})$ and then over the fiber of $L$.
By equation (\ref{eq:useful-trans}),
the volume form of $(S^{m})^{k+1}$ at $A$ can be written as
\[
d(S^{m})^{k+1} = \operatorname{\mathsf{vol}}(A)^{m-k} \,d\mathrm{G}_{k+1}(\mathbb R^{m+1})(L)
\wedge dS(L)^{k+1},
\]
where $S(L)^{k+1}$ denotes $(k+1)$-fold product of the unit sphere of
$L$. By invariance under orthogonal transformations, the integral over
the fiber does not depend on~$L$. We may therefore assume that
$L=\mathbb R^{k+1}$, in which case the fiber over $L$ can be identified
with $\mathcal {R}_{k+1,k}^k(\varepsilon)$. Thus we conclude from
equation (\ref{eq:weight}) that
%
%
\begin{equation}\label{eq:nice}
\frac{\operatorname{\mathsf{vol}}\mathcal{R}_{n,m}^{k}(\varepsilon)}
{\mathcal O_m^{n-k-1}} = \mathcal G_{k+1,m+1}
\int_{A\in\mathcal{R}_{k+1,k}^k(\varepsilon)}
\operatorname{\mathsf{vol}}(A)^{m-k}
\lambda_m(t(A))^{n-k-1}\,d(S^{k})^{k+1}.\hspace*{-29pt}
\end{equation}

In a next step, we will perform a change of variables in order to express
the integral on the right-hand side of equation (\ref{eq:nice}) as an
integral over $t$
involving the coefficients $C(m,k)$.

Note that by Lemma \ref{le:pt}, $p(A)\in S^k$ and $t(A)\in(0,1)$ are
defined for any $A\in\mathrm{GL}(k+1,\mathbb R)$ and depend smoothly on
$A$. A moment's thought (together with Lemmas \ref{le:crit} and
\ref{le:elusive}) shows that we have the following complete
characterization of $\mathcal{R}_{k+1,k}^k(\varepsilon)$:
\begin{eqnarray*}
\mathcal{R}_{k+1,k}^k(\varepsilon) &=& \{ A \in(S^k)^{k+1} \mid
\mbox{$A$ is centered with respect to $t(A)p(A)$,} \\
&&\hspace*{86pt} 0 < t(A) <\varepsilon, \forall i\ \langle a_i,p(A)\rangle= t(A) \}.
\end{eqnarray*}
For $A\in\mathcal{R}_{k+1,k}^k(\varepsilon)$ set $p:=p(A)$, $t:=t(A)$,
and $r := r(t):=(1-t^2)^{1/2}$. For $i\in[k+1]$ we define $b_i$ as the
scaled-to-one orthogonal projection of $a_i$ onto the orthogonal
complement of $\mathbb Rp$, briefly $a_i=rb_i +tp$. The matrix $B=B(A)$
with the rows $b_1,\ldots,b_{k+1}$ can be written as $B =\frac{1}{r}
(A-tep^\top)$. Clearly, $B$ is centered.

We define now
\[
W_{k}= \{(B,p)\in(S^{k})^{k+1}\times S^{k} \mid Bp=0
\mbox{ and } B \mbox{ is centered} \}.
\]
This is a Riemannian submanifold of $(S^{k})^{k+2}$ of
dimension $k(k+1)-1$. We thus have a map,
\[
\varphi_k\dvtx\mathcal{R}_{k+1,k}^k(\varepsilon) \rightarrow W_k
\times(0,\varepsilon),\qquad
A \mapsto(B(a),p(A),t(A)).
\]
The inverse of this map is given by $(B,p,t)\mapsto A=rB+ tep^\top$. It
is well defined since, by Lemma \ref{le:elusive}, $A$ is invertible
when $B$ is centered. The Jacobian $J(A)$ of the diffeomorphism
$\varphi_k$ is stated in the next lemma, whose proof will be
momentarily postponed. We remark that this lemma can also be derived
from \cite{reitz}, Lemma 1 (a~special case of Z\"{a}hle's theorem
\cite{zaehle}) with $K$ being the unit ball.
\begin{lemma}\label{le:Jacobian}
The volume form of $(S^{k})^{k+1}$ at $A$ can be expressed in terms of
the volume form of $W_k\times(0,\varepsilon)$ as follows:
\[
d(S^{k})^{k+1} = J(A) \,dW_k\wedge dt = \frac{r^{(k-2)(k+1)}
\operatorname{\mathsf{vol}}(A)}{t} \,dW_k\wedge dt .
\]
\end{lemma}

We now express the Jacobian $J(A)$ in terms of $(B,p,t)$. The volume
$\operatorname{\mathsf{vol}}(A)$ of the parallelepiped spanned by the
$a_i$ equals $(k+1)!$ times the volume of the pyramid with apex $0$ and
base $\Delta(a_1,\ldots,a_{k+1})$, the latter denoting the simplex with
vertices $a_1,\ldots,a_{k+1}$. Moreover, this pyramid has height $t$
and it is well known that the volume of a $(k+1)$-dimensional pyramid
with height $t$ and base $B$ is $\frac{t}{k+1}$ times the
($k$-dimensional) volume of $B$. This implies
\[
\operatorname{\mathsf{vol}}(A) = (k+1)! \frac{t}{k+1}
\operatorname{\mathsf{vol}}_{k-1}\Delta (rb_1,\ldots,rb_{k+1}) = k! r^k
t \operatorname{\mathsf{vol}}_{k}\Delta(B).
\]
From this expression, together with (\ref{eq:nice}), we conclude that
\begin{eqnarray*}
\frac{\operatorname{\mathsf{vol}}\mathcal{R}_{n,m}^{k}(\varepsilon
)}{\mathcal O_m^{n}} &=& \frac{\mathcal G _{k+1,m+1}}{\mathcal O_m^{k+1}}
\int_{W_k\times(0,\varepsilon)} \operatorname{\mathsf{vol}}(A)^{m-k+1}
\frac {r(t)^{(k-2)(k+1)}}{t}\\
&&\hspace*{84.6pt}{}\times
\lambda_m(t)^{n-k-1} \,d(S^k)^{k+1} \\
&=& \frac{\mathcal G_{k+1,m+1}(k!)^{m-k+1}}{\mathcal O_m^{k+1}}\int_{W_k}
\operatorname{\mathsf{vol}}_k\Delta(B)^{m-k+1} \,dW_k\\
&&{}\times \int_{0}^{\varepsilon}
t^{m-k}r(t)^{km-2}\lambda_m(t)^{n-k-1} \,dt.
\end{eqnarray*}
Consider the projection $\mathrm{pr}\dvtx W_k\to S^k, (B,p)\mapsto p$.
We note that its fiber over~$p$ can be identified with the set $M_k$
(cf. Definition \ref{def:C}). By the invariance of
$\operatorname{\mathsf{vol}}_k \Delta(B)$ under rotation of $p\in S^k$,
we get
\[
\int_{W_k} \operatorname{\mathsf{vol}}_k\Delta(B)^{m-k+1}\, dW_k =
\mathcal O_k \int_{M_k}\operatorname{\mathsf{vol}}_k\Delta(B)^{m-k+1}\,
d(S^{k-1})^{k+1}.
\]
Note that, up to a scaling factor, the right-hand side above is the
coefficient $C(m,k)$ introduced in Definition \ref{def:C}. Using
$(\mathcal O_k/\mathcal O_m)\mathcal G_{k+1,m+1}=\mathcal G_{k,m}$ we
obtain
%
%
\begin{equation}\label{eq:reuse}
\frac{\operatorname{\mathsf{vol}}\mathcal{R}_{n,m}^{k}(\varepsilon)}{\mathcal
O_m^{n}} = C(m,k) \int_{0}^{\varepsilon} t^{m-k}r(t)^{km-2}\lambda
_m(t)^{n-k-1}\, dt.
\end{equation}
This completes the proof.
\end{pf}
\begin{pf*}{Proof of Lemma \protect\ref{le:Jacobian}}
For given $p\in S^k$, let $S(p^{\perp})$ denote the $(k-1)$-subsphere
of $S^k$ perpendicular to $p$. At a point $(B,p)\in W_k$ we have
$dW_k=dS(p^{\perp})^{k+1}\wedge dS^k$, and we have
$dS(p^{\perp})^{k+1}=dS(p^{\perp})\wedge\cdots\wedge dS(p^{\perp})$ at
the point $B=(b_1,\ldots,b_{k+1})$. The Jacobian $J(A)$ we are looking
for is hence determined by
\[
d(S^{k})^{k+1}(A) = J(A) \,dS(p^{\perp})^{k+1}(B) \wedge dS^k (p) \wedge dt.
\]

Choose an orthonormal moving frame $e_1,\ldots,e_k,e_{k+1}$ with
$p(A)=e_{k+1}$. For $1\leq i\leq k$ define the one-forms
$\mu_{i}:=-\langle e_i,dp\rangle$ (compare Section \ref{se:Grass}).
Then the volume form of $S^k$ at $p$ is given by
$dS^k(p)=\mu_1\wedge\cdots\wedge \mu_k$.

Differentiating $te=Ap$ we get $e \,dt=dA p+A\, dp$. By multiplying both
sides with $A^{-1}$ and using formula (\ref{eq:pt-formula}) we obtain
\[
p (dt/t) - dp = A^{-1}\,dA p.
\]
Let $Q$ denote the $(k+1)\times(k+1)$ matrix having the $e_i$ as rows.
With respect to this basis, the above equation takes the form
%
%
\begin{equation}\label{eq:1matrix-D}
(\mu_1,\ldots,\mu_k,d t/t)^{\top} = Q \bigl(p (dt/t)-dp \bigr) = Q A^{-1} \,dA p.
\end{equation}
Wedging the entries on both sides yields
%
%
\begin{equation}\label{eq:wedge1}dS^k(p)\wedge dt
= t \operatorname{\mathsf{vol}}(A)^{-1}\langle
p,da_1\rangle\wedge\cdots\wedge \langle p,da_{k+1}\rangle.
\end{equation}
The volume form of $S(p^{\perp})$ at $b_i$
is given by
\[
dS(p^{\perp})=\langle e_1,db_i\rangle\wedge\cdots\wedge\langle
e_{k-1},db_i\rangle.
\]
In order to compare $dS(p^{\perp})^{k+1}\wedge dS^k\wedge dt$ with
$d(S^k)^{k+1}$ we use a different moving frame.
Fix an $i$, $1\leq i\leq k+1$, and choose the moving frame as above
and additionally with $e_k=b_i$.
Consider the modified frame
$\tilde{e}_1,\ldots,\tilde{e}_{k+1}$ that arises after rotating
$b_i$ to $a_i$
and leaving the orthogonal complement of
$\mathrm{span}\langle a_i,p\rangle$ fixed, that is, $\tilde
{e}_j=e_j$ for $1\leq
j\leq k-1$, $\tilde{e}_k:=a_i$,
and $\tilde{e}_{k+1}=-t b_i + r p$
(cf. Figure \ref{fig:2}).

%
%
\begin{figure}[b]

\includegraphics{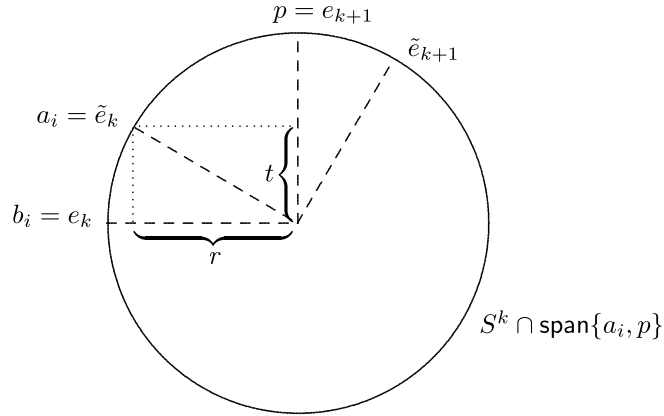}

\caption{The frame $(e_i)$ and its modification $(\tilde
{e}_i)$.}\label{fig:2}
\end{figure}

This implies
$\langle\tilde{e}_{k+1},da_i\rangle=r\langle p,da_i\rangle-t\langle
b_i,da_i\rangle=(1/r) \langle p,da_i\rangle$
where we have used that $b_i=(a_i-tp)/r$ for the last equality.
Hence the volume form of $S^k$ at $a_i$ equals
\[
dS^k(a_i)=(1/r)\langle e_1,da_i\rangle\wedge\cdots\wedge
\langle e_{k-1},da_i\rangle\wedge\langle p,da_i\rangle.
\]
If we wedge the $\langle e_1,da_i\rangle\wedge\cdots\wedge\langle
e_{k-1},da_i\rangle$ to both sides of equation (\ref{eq:wedge1}),
we obtain, on the right-hand side,
\[
\frac{t}{\operatorname{\mathsf{vol}}(A)}\bigwedge_{i=1}^{k+1}
\langle e_1,da_i\rangle\wedge\cdots\wedge\langle
e_{k-1},da_i\rangle\wedge\langle p,da_i\rangle
= \frac{t}{\operatorname{\mathsf{vol}}(A)} r^{k+1}\, d(S^k)^{k+1}(A).
\]
On the left-hand side we obtain, using
$\langle e_j,da_i\rangle=r\langle e_j,db_i\rangle+t\langle
e_j,dp\rangle$
and taking into account that
$\langle e_j,dp\rangle\wedge dS^k(p)=0$ since $dS^k(p)=\bigwedge
_{j=1}^k\langle e_j,dp\rangle$,
\begin{eqnarray*}
\bigwedge_{i=1}^{k+1}\bigwedge_{j=1}^{k-1} \langle e_j,da_i\rangle
\wedge dS^k\wedge dt
&=& r^{(k-1)(k+1)}\bigwedge_{i=1}^{k+1}\bigwedge_{j=1}^{k-1}
\langle e_j,db_i\rangle\wedge dS^k\wedge dt\\
&=& r^{(k-1)(k+1)} dS(p^\perp)^{k+1}\wedge dS^k\wedge dt.
\end{eqnarray*}
This implies that
$J(A) = t^{-1} r^{(k-2)(k+1)}\operatorname{\mathsf{vol}}(A)$ as claimed.
\end{pf*}

\subsection{The infeasible case}

Recall that $\mathcal{I}_{n,m}= (S^{m})^n \setminus\mathcal
{F}_{n,m}$ denotes
the set of infeasible instances. We define, for $I\subseteq[n]$,
\[
\mathcal{I}_{n,m}^I(\varepsilon) := \{A\in\mathcal{I}_{n,m} \mid
\mbox{${\mathscr C}(A)> \varepsilon^{-1}$ and $A$ has a SIC with
blocking set $I$} \}.
\]
We note that by symmetry, the volume of $ \mathcal
{I}_{n,m}^I(\varepsilon)$ only depends
on the cardinality of $I$.
\begin{lemma}\label{lem:if_sic}
$\mathcal{I}_{n,m}^I(\varepsilon)$ has measure zero if $|I|>m+1$.
\end{lemma}
\begin{pf}
If $A\in\mathcal{I}_{n,m}^I(\varepsilon)$, then $\{a_i\mid i\in I\}$
is contained in the boundary of a SIC of $A$ with blocking set $I$.
Hence the affine hull of $(a_i)_{i\in I}$ has dimension less than~$m$.
However, if $|I|>m+1$, the latter dimension is almost surely $m+1$.
\end{pf}

It is known \cite{ChCH05}, Proposition 4.2, that in
the infeasible case, blocking sets have at least $m+1$ elements.
This fact, together with Lemma \ref{lem:if_sic}, implies that
%
%
\begin{equation}\label{eq:inf}
\operatorname{\mathsf{vol}}\mathcal{I}_{n,m}(\varepsilon) \le\pmatrix{n\cr
{m+1}} \operatorname{\mathsf{vol}}\mathcal{I}
_{n,m}^{[m+1]}(\varepsilon).
\end{equation}
As with $\mathcal{F}_{n,m}$, for ease of notation, we write
$\mathcal{I}_{n,m}^m(\varepsilon):=\mathcal
{I}_{n,m}^{[m+1]}(\varepsilon)$.

The inequality in Theorem \ref{th:main_lp} for the infeasible case
follows immediately from (\ref{eq:inf}) and the following proposition.
\begin{proposition}\label{th:main-if}
We have for $\varepsilon\in(0,1]$,
\[
\frac{\operatorname{\mathsf{vol}}\mathcal{I}_{n,m}^m(\varepsilon)}{\mathcal O_m^n}
\leq C(m,m) \int_0^{\varepsilon} (1-t^2)^{({m^2-2})/{2}}
\bigl(1-\lambda_m(t) \bigr)^{n-m-1} \,dt.
\]
\end{proposition}
\begin{pf}
Consider the projection
\[
\psi\dvtx\mathcal{I}_{n,m}^m(\varepsilon) \to(S^m)^{m+1},\qquad
A'=(a_1,\ldots,a_n)\mapsto A=(a_1,\ldots,a_{m+1}) .
\]
To investigate the image and the fibers of $\psi$ assume
$A'\in\mathcal{I}_{n,m}^m(\varepsilon)$. Then there exists $p\in S^m$
and $\alpha\in(\pi/2,\pi]$ such that $\mathsf{cap}(p,\alpha)$ is a SIC
of $A'\in\mathcal{I}_{n,m}^m(\varepsilon)$ with blocking set $[m+1]$.
Then we have that $\langle a_i,p\rangle= t$ for all $i\in[m+1]$ and
$\langle a_i,p\rangle> t$ for all $i\in[n]\setminus[m+1]$ where
$t:=\cos\alpha\in[-1,0)$. Lemma \ref{lem:new} implies that
$tp\in\mathsf{conv}\{a_1,\ldots,a_{m+1}\}$. In turn, Lemma
\ref{le:crit} implies that $\mathsf{cap}(-p,\pi-\alpha)$ is a SIC of
$A$ with blocking set $[m+1]$, and we obtain that
$A\in\mathcal{F}_{m+1,m}^m(\varepsilon)$. These reasonings show that
the image of $\psi$ is contained in
$\mathcal{F}_{m+1,m}^m(\varepsilon)$.

Suppose now that $a_1,\ldots,a_{m+1}$ are linearly independent. Then it
follows from Lemma \ref{le:pt} that the vector $p$ is uniquely
determined by $A$. This implies that the fiber of $A$ under $\psi$ is
contained in $\{A\}\times\mathsf{cap}(p,\alpha)^{n-m-1}$. We conclude
that for almost all $A\in\mathcal{F}_{m+1,m}^m$
\[
\frac{\operatorname{\mathsf{vol}}\psi^{-1}(A)}{\mathcal O_m^{n-m-1}} \le
\bigl(1-\lambda_m(t)\bigr)^{n-m-1} .
\]
From these observations we obtain, by Fubini,
\[
\frac{\operatorname{\mathsf{vol}}\mathcal{I}_{n,m}^{m}(\varepsilon)}{\mathcal
O_{m}^{n-m-1}} \leq\int_{A\in\mathcal{F}_{m+1,m}^m(\varepsilon)}
\bigl(1-\lambda_m(t(A))\bigr)^{n-m-1}\, d(S^m)^{m+1} .
\]
In the proof of Proposition \ref{th:main-f}
we derived, from the integral representation (\ref{eq:weight})
for $\frac{\operatorname{\mathsf{vol}}\mathcal{F}_{n,m}^{k}(\varepsilon
)}{\mathcal O_m^{n-k-1}}$,
formula (\ref{eq:reuse}).
In exactly the same way we can show that
\[
\frac{\operatorname{\mathsf{vol}}\mathcal{I}_{n,m}^{m}(\varepsilon)}{\mathcal
O_m^{n}}\leq
C(m,m) \int_0^{\varepsilon} (1-t^2)^{({m^2-2})/{2}}
\bigl(1-\lambda_m(t) \bigr)^{n-m-1} \,dt,
\]
which proves the assertion.
\end{pf}
\begin{remark}\label{rem:stevens}
(i) It may be of interest to compare, in the case $m=1$, the upper bound
for $p(n,1,\alpha)$ which follows from our results with
the exact expression (\ref{eq:stevens})
for this quantity shown by Stevens.
Recall that the latter is
\begin{eqnarray*}
p(n,1,\alpha) &=& n \biggl(1-\frac{\alpha}{\pi} \biggr)^{n-1}-
\pmatrix{n\cr2} \biggl(1-\frac{2\alpha}{\pi} \biggr)^{n-1}+ \cdots\\
&&{} + (-1)^{k+1}\pmatrix{n\cr k} \biggl(1-\frac{k\alpha}{\pi} \biggr)^{n-1},
\end{eqnarray*}
where $k=\lfloor\frac{\pi}{\alpha}\rfloor$.
For $\alpha\in[0,\pi/2]$, Propositions \ref{prop:equiv}
and \ref{th:main-if} yield
\begin{eqnarray*}
p(n,1,\alpha)&=& \mathsf{Prob}\{A\in\mathcal{F}_{n,1}\}
+\mathsf{Prob}\biggl\{A\notin\mathcal{F}_{n,1}
\mbox{ and }
{\mathscr C}(A)\geq\frac{1}{\cos(\alpha)} \biggr\}\\
&\leq& \frac{n}{2^{n-1}}+
\pmatrix{n\cr2}C(1,1)\int_0^{\cos\alpha}
(1-t^2)^{-1/2} \bigl(1-\lambda_1(t) \bigr)^{n-2} \,dt.
\end{eqnarray*}
We use now that $C(1,1)=\frac2{\pi}$, as
shown in Example \ref{ex:k=2}. Then
\begin{eqnarray*}
p(n,1,\alpha)
& =& \frac{n}{2^{n-1}}+
\frac{n(n-1)}{\pi}
\int_0^{\cos\alpha} (1-t^2)^{-1/2}
\biggl(1-\frac{\arccos t}{\pi} \biggr)^{n-2}\, dt\\
& \le& \frac{n}{2^{n-1}}+
n(n-1)\int_{\alpha/\pi}^{1/2}
(1-x )^{n-2}\,dx\\
&=& n \biggl(1-\frac{\alpha}{\pi} \biggr)^{n-1}.
\end{eqnarray*}
That is, we get Stevens's first term.\smallskipamount=0pt
\begin{longlist}[(ii)]
\item[(ii)]
It may also be of interest to compare, for the case $m=2$, our upper bound
for $p(n,2,\alpha)$ with the upper bound in (\ref{eq:Gilbert}) obtained
by Gilbert \cite{Gilbert65}. Recall that the latter gives
\[
p(n,2,\alpha)\leq\tfrac43n(n-1)\lambda(1-\lambda)^{n-1},
\]
where $\lambda$ denotes the fraction of the surface of the sphere
covered by the
cap of radius $\alpha$. It is easy to see that our bound implies
\[
p(n,2,\alpha)\leq\frac{1}{2^n}(n^2-n+2)+\frac12n(n-2)(1-\lambda)^{n-1}.
\]
The first term in this sum is negligible for large $n$. The second term
compares with
Gilbert's for moderately large caps but it becomes considerably worse
for small
caps.
\end{longlist}
\end{remark}

\section{On the values of the coefficients $C(m,k)$}\label{se:C(m,k)}

In this section we provide estimates for the numbers $C(m,k)$. In
Section \ref{sec:Cmk_bounds} we derive upper and lower bounds for them
which are elementary functions in $m$ and $k$. In the case $m=k$ the
upper bound is actually an equality, yielding an exact expression for
$C(m,m)$. Then in Section \ref{sec:proof_corol} we use these bounds to
prove Theorems \ref{cor:2} and \ref{corol:1}. Finally, in Section
\ref{sec:Cmk_part} we briefly describe how to derive an exact
expression for $C(m,m-1)$ and how, for any given $m$, one may obtain
the values of the $C(m,k)$, $k=1,\ldots,m$, by solving an $m\times m$
linear system.

\subsection{Bounding the coefficients $C(m,k)$}\label{sec:Cmk_bounds}

Our first result provides bounds for $C(m,k)$ in terms of
volumes of spheres.
\begin{lemma}\label{le:bound}
We have for $1\le k\le m$,
\[
\frac{k+1}{2^k}\frac{\mathcal O_{k-1}\mathcal O_{m-k}}{\mathcal O_{m}}
\leq C(m,k)\leq \frac{(k+1)^{m-k+1}}{2^k}\frac{\mathcal O_{k-1}\mathcal
O_{m-k}}{\mathcal O_{m}}
\]
with equalities if $k=m$. In particular,
$C(m,m)=\frac{m+1}{2^{m-1}}\frac{\mathcal O_{m-1}}{\mathcal O_m}$.
\end{lemma}
\begin{pf}
Recall Definition \ref{def:C} of the $C(m,k)$,
\[
C(m,k)=\frac{(k!)^{m-k+1}}{\mathcal O_m^k} \mathcal G_{k,m} \int_{M_k}
(\operatorname{\mathsf{vol}}\Delta)^{m-k+1} \,dS .
\]
We set $S:=(S^{k-1})^{k+1}$ and denote by $U$ the open dense subset of
$S$ consisting of all $B=(b_1,\ldots,b_{k+1})$ such that every $k$ of
these vectors are linearly independent. By Definition
\ref{def:centered}, $M_k$ is contained in $U$.

Set $\Delta(B)=\mathsf{conv}\{b_1,\ldots,b_{k+1}\}$ and
$\Delta_i(B)=\mathsf{conv}(0,b_1,\ldots,\hat{b}_i,\ldots,b_{k+1})$\break
(where $\hat{b}_i$ means that $b_i$ is omitted). We define, for $B\in
U$,
\[
\mathsf{mvol}(B):=\sum_{i=1}^{k+1} \operatorname{\mathsf{vol}}\Delta_i(B).
\]
For $B\in M_k$ we clearly have
$\mathsf{mvol}(B)=\operatorname{\mathsf{vol}}\Delta (B)$, but in
general this is not the case.

The essential observation is now the following:
%
%
\begin{equation}\label{eq:trick}
\int_{M_k}(\operatorname{\mathsf{vol}}\Delta)^{m-k+1}\,dS =
\frac{1}{2^k}\int_{U}\mathsf{mvol}^{m-k+1}\,dS.
\end{equation}
In order to show this, note that for $B\in U$ there exists a unique
$\mu\in\mathbb R^{k+1}$ with $\mu_{k+1}=1$, $\mu_1\cdots\mu_k\ne0$, and
such that $\sum_{i=1}^{k+1}\mu_i b_i=0$. This allows to define the map
$\phi\dvtx U\rightarrow\{-1,1\}^k,
B\mapsto(\mathsf{sgn}(\mu_1),\ldots,\mathsf{sgn}(\mu_{k}))$. Note that
$M_k=\phi^{-1}(1,\ldots,1)$. Moreover, each $\sigma\in\{-1,1\}^k$
defines an isometry,
\[
M_k\rightarrow\phi^{-1}(\sigma),\qquad B\mapsto
\sigma B :=(\sigma_1b_1,\ldots,\sigma_kb_k,b_{k+1}).
\]
It follows that $\mathsf{mvol}(B)=\mathsf{mvol}(\sigma B)$ since
changing the signs of rows does not alter the absolute values of
determinants. This implies
\begin{eqnarray*}
\int_{U}\mathsf{mvol}^{m-k+1} \,dS &=& \sum_{\sigma\in\{-1,1\}^k}
\int_{\phi^{-1}(\sigma)} \mathsf{mvol}^{m-k+1} \,dS\\
&=& 2^k \int_{M_k} \mathsf{mvol}^{m-k+1} \,dS,
\end{eqnarray*}
which proves the claimed equation (\ref{eq:trick}).

Recall now the norm inequalities
%
%
\begin{eqnarray}\label{eq:ineq}
(x_1^\ell+\cdots+x_p^\ell) &\leq&
(x_1+\cdots+x_p)^\ell\nonumber\\[-8pt]\\[-8pt]
&\leq& p^{\ell-1} (x_1^\ell+\cdots+x_p^\ell)\qquad
\mbox{for $x_i\ge0, \ell\ge1$},\nonumber
\end{eqnarray}
where the last follows from the convexity of the function $\mathbb
R\to\mathbb R, y\mapsto y^\ell$. For the upper bound in the statement
we now estimate the right-hand side of equation (\ref{eq:trick}) using the last
inequality above (with $p=k+1$ and $\ell=m-k+1$). We obtain
\begin{eqnarray*}
\int_{S} \mathsf{mvol}^{m-k+1}\,dS
&\leq& (k+1)^{m-k}\sum_{i=1}^{k+1}\int_{S}(\operatorname{\mathsf{vol}}\Delta_i)^{m-k+1}
\,dS\\
&=& (k+1)^{m-k+1}\int_{S}(\operatorname{\mathsf{vol}}\Delta_{k+1})^{m-k+1} \,dS\\
&=& \frac{(k+1)^{m-k+1}}{k!^{m-k+1}}
\int_{S} |{\det\tilde{B}}|^{m-k+1} \,dS,
\end{eqnarray*}
where $\tilde{B}\in\mathbb R^{k\times k}$ denotes the matrix with rows
$b_1,\ldots,b_k$. Since the integrand on the right does not depend on
$b_{k+1}$, we can integrate over $\tilde{B}\in(S^{k-1})^k$ and pull out
a factor $\mathcal O_{k-1}$ obtaining
\begin{eqnarray*}
\int_{S} \mathsf{mvol}^{m-k+1}\,dS
&\leq& \frac{(k+1)^{m-k+1}}{k!^{m-k+1}} \mathcal O_{k-1} \int_{(S^{k-1})^{k}}
|{\det\tilde{B}}|^{m-k+1} \,d(S^{k-1})^{k}\\
&=& \frac{(k+1)^{m-k+1}}{k!^{m-k+1}} \mathcal O_{k-1}^{k+1}
\mathbf{E}(|{\det\tilde{B}}|^{m-k+1} ) .
\end{eqnarray*}
We plug in here the formula of the moments from Lemma \ref{le:highmom}.
Putting everything together, and using $\mathcal G_{k,m}=(\mathcal
O_{m-k}/\mathcal O_{m})\mathcal G_{k,m+1}$, the claimed upper bound on
$C(m,k)$ follows. The lower bound is obtained by doing the same
reasoning but now using the left-hand side inequality in
(\ref{eq:ineq}).

In the case $k=m$ upper and lower bounds coincide and we get
equalities for $C(m,m)$.
\end{pf}
\begin{remark}
In the case $k=1$ the upper bound
in Lemma \ref{le:bound} coincides with the value for $C(m,1)$
shown in Example \ref{ex:k=2}.
\end{remark}

For the proofs of Theorems \ref{cor:2} and \ref{corol:1}
we need a more explicit expression for the bounds on the $C(m,k)$.
We devote the rest of this section to deriving
such expressions.
\begin{lemma}\label{le:easy}
For $1\leq k \leq m$ we have
\[
\frac{\mathcal O_{k-1}\mathcal O_{m-k}}{\mathcal O_{m}} \leq\sqrt
{\frac{\pi}{2}}
k^{3/4} \sqrt{\pmatrix{m\cr k}} .
\]
In the cases $k=1$ or $k=m$ one has the sharper bound
$\frac{2\mathcal O_{m-1}}{\mathcal O_m} \leq\sqrt{m}$.
\end{lemma}

The proof uses bounds on Gamma functions, which we derive next.
\begin{lemma}\label{lem:gamma}
For all $r\geq1$,
\begin{eqnarray*}
r^{1/4} 2^{-({r-1})/{2}} \sqrt{(r-1)!}
&\leq&\Gamma\biggl(\frac{r+1}{2} \biggr)\\
&\leq&\sqrt{\frac{\pi}{2}} r^{1/4} 2^{-({r-1})/{2}}
\sqrt{(r-1)!} .
\end{eqnarray*}
\end{lemma}
\begin{pf}
The double factorials $k!!$ are defined as follows.
For $k$ even, $k!! :=k(k-2)(k-4) \cdots2$,
and for $k$ odd, $k!! :=k(k-2)(k-4) \cdots3\cdot1$. Also,
by convention, $0!!=1$.
By the functional equation $\Gamma(x+1)=x\Gamma(x)$
of the Gamma function,
it easily follows that for $r\in\mathbb N$, $r\geq1$,
%
%
\begin{equation}\label{eq:gamma}
\Gamma\biggl(\frac{r+1}{2} \biggr) =
\cases{
\sqrt{\dfrac{\pi}{2}} (r-1)!! 2^{-({r-1})/{2}},
&\quad if $r$ is even,\cr
(r-1)!! 2^{-({r-1})/{2}},
&\quad if $r$ is odd.}
\end{equation}
We estimate now double factorials in terms of factorials.
If $r\geq2$ is even,
\begin{eqnarray*}
(r!!)^2 &=& r r (r-2) (r-2) \cdots4\cdot4\cdot2\cdot2\\
&=& r(r-1)\frac{r}{r-1} (r-2)(r-3)\frac{r-2}{r-3}\cdots4\cdot3
\frac43 2\cdot2\\
&=& r! \frac{r}{r-1} \frac{r-2}{r-3}\cdots\frac43 2\\
&=& r! \sqrt{\frac{r}{r-1} \frac{r}{r-1} \frac{r-2}{r-3} \frac{r-2}{r-3}
\cdots\frac43 \frac43 2\cdot2}.
\end{eqnarray*}
We use that
$\frac{\ell+1}{\ell}\leq\frac{\ell}{\ell-1}$
for $\ell\geq2$ to deduce from this
%
%
\begin{equation}\label{eq:df1}
r!\sqrt{r+1}\leq(r!!)^2 \leq r!\sqrt{2r}\qquad
\mbox{for $r\geq2$ even.}
\end{equation}
Similarly, for $r\geq1$ odd, one shows that
\[
(r!!)^2=r! \sqrt{\frac{r}{r-1} \frac{r}{r-1} \frac{r-2}{r-3} \frac
{r-2}{r-3}
\cdots\frac54 \frac54 \frac32 \frac32},
\]
which implies
%
%
\begin{equation}\label{eq:df2}
r! \sqrt{\frac{r+1}{2}}\leq(r!!)^2 \leq r!\sqrt{r}\qquad
\mbox{for $r\geq1$ odd.}
\end{equation}
By applying the bounds (\ref{eq:df1}) and (\ref{eq:df2}) to (\ref{eq:gamma})
and noting that $2^{1/4}\le\sqrt{\frac{\pi}{2}}$, the claim follows.
\end{pf}
\begin{pf*}{Proof of Lemma \protect\ref{le:easy}}
Assume that $2\leq k <m$. Then, using Lemma \ref{lem:gamma}, we deduce that
\begin{eqnarray*}
\frac{\mathcal O_{k-1}\mathcal O_{m-k}}{\mathcal O_{m}}
&=&
2\frac{\Gamma(({m+1})/{2} )}
{\Gamma({k}/{2} )\Gamma(({m-k+1})/{2} )}\\
&\leq& \sqrt{\frac{\pi}{2}}
\sqrt{\frac{(m-1)!}{(k-2)!(m-k-1)!}}
\biggl(\frac{m}{(k-1)(m-k)} \biggr)^{1/4}\\
&= & \sqrt{\frac{\pi}{2}}
\sqrt{\pmatrix{m\cr k}}
\sqrt{\frac{(m-k)k(k-1)}{m}}
\biggl(\frac{m}{(k-1)(m-k)} \biggr)^{1/4}\\
&\leq& \sqrt{\frac{\pi}{2}}
\sqrt{\pmatrix{m\cr k}} k^{3/4}.
\end{eqnarray*}
The cases $k=1$ and $k=m$ follow similarly from Lemma \ref{lem:gamma}.
\end{pf*}

An immediate consequence of Lemmas \ref{le:easy} and \ref{lem:gamma}
are the following bounds on the coefficients $C(m,k)$.
\begin{proposition}\label{prop:conj}
For $1\le k<m$,
\[
C(m,k)\le\sqrt{\frac{\pi}{2}}
\frac{(k+1)^{m-k+1}}{2^k} k^{3/4} \sqrt{\pmatrix{m\cr k}}.
\]
In addition, for all $m\geq1$,
\[
C(m,m)\leq\frac{(m+1)\sqrt{m}}{2^m}.
\]
\end{proposition}
\begin{remark}
Using Lemmas \ref{le:bound} and \ref{lem:gamma} it is easy
to obtain lower bounds for the $C(m,k)$ similar to the upper bounds
in Proposition \ref{prop:conj}.
\end{remark}

\subsection[Proof of Theorems 1.2 and 1.4]{Proof of
Theorems \protect\ref{cor:2} and \protect\ref{corol:1}}
\label{sec:proof_corol}

The following identity is repeatedly used in the proof.
\begin{lemma}\label{lem:series}
We have
$\sum_{n=k}^{\infty}{n\choose k} z^{n-k}=(1-z)^{-k-1}$
for $k\in\mathbb N$ and $z\in(0,1)$.
\end{lemma}
\begin{pf}
Take the $k$th derivative on both sides of
$\sum_{n=0}^\infty z^n = \frac{1}{1-z}$.
\end{pf}
\begin{pf*}{Proof of Theorem \protect\ref{cor:2}}
By definition, we have $N(m,\alpha)>n$ iff
$\mathsf{cap}(a_1,\break\alpha)\cup\cdots\cup\mathsf{cap}(a_n,\alpha
)\ne S^n$.
Hence
\[
\mathbf{E}(N(m,\alpha)) = \sum_{n=0}^\infty\mathsf
{Prob}\bigl( N(m,\alpha) > n\bigr)
= \sum_{n=0}^\infty p(n,m,\alpha) .
\]
We assume that $\alpha\le\pi/2$.
Since $p(n,m,\alpha)=1$ for $n\le m+1$, we conclude
%
%
\begin{equation}\label{eq:stieljes}
\mathbf{E}(N(m,\alpha))= m +1 + \sum_{n=m+1}^\infty p(n,m,\alpha).
\end{equation}
Proposition \ref{prop:equiv} states that, for $\alpha\in(0,\frac
{\pi}{2}]$,
and $\varepsilon=\cos(\alpha)$
\[
p(n,m,\alpha)= 2^{1-n} \sum_{k=0}^{m}\pmatrix{n-1\cr k} +
P_{n,m}(\varepsilon),
\]
where we have put
\[
P_{n,m}(\varepsilon) :=
\mathsf{Prob}\{A\in\mathcal{I}_{n,m}\mbox{ and }
{\mathscr C}(A)\ge\varepsilon^{-1}\}.
\]

We first estimate
\[
T:= \sum_{n=m+1}^\infty2^{1-n} \sum_{k=0}^{m}\pmatrix{n-1\cr k}
\]
as follows (take $r=n-1$)
\[
T = \sum_{k=0}^{m}\sum_{r=m}^\infty\pmatrix{r\cr k}
\biggl(\frac12 \biggr)^r
\leq\sum_{k=0}^{m} \biggl(\frac12 \biggr)^k
\sum_{r=k}^\infty\pmatrix{r\cr k} \biggl(\frac12 \biggr)^{r-k}
- \sum_{k=0}^{m-1} \frac1{2^k}.
\]
Applying Lemma \ref{lem:series} to the last expression we obtain
\[
T \le\sum_{k=0}^{m} \biggl(\frac12 \biggr)^k 2^{k+1}
- 2 +\frac1{2^{m-1}} \le2m +1.
\]

We now estimate
$T^* := \sum_{n=m+1}^\infty P_{n,m}(\varepsilon)$
using Theorem \ref{th:main_lp} which tells us that
\[
P_{n,m}(\varepsilon) \le\pmatrix{n\cr{m+1}} C(m,m)
\int_0^\varepsilon(1-t^2)^{({m^2-2})/{2}}\bigl(1-\lambda
_m(t)\bigr)^{n-m-1} \, dt.
\]
Hence, using Lemma \ref{lem:series} again,
\begin{eqnarray*}
T^* &=& C(m,m)\sum_{n=m+1}^\infty
\pmatrix{n\cr{m+1}} \int_0^\varepsilon(1-t^2)^{({m^2-2})/{2}}
\bigl(1-\lambda_m(t)\bigr)^{n-m-1} \,dt\\
&\leq& C(m,m)\int_0^\varepsilon\sum_{n=m+1}^\infty
\pmatrix{n\cr{m+1}} \bigl(1-\lambda_m(t)\bigr)^{n-m-1} \,dt\\
& = & C(m,m)\int_0^\varepsilon\lambda_m(t)^{-m-2} \,dt
\leq C(m,m) \frac{\varepsilon}{\lambda(\varepsilon)^{m+2}} .
\end{eqnarray*}
Plugging in the estimate for $C(m,m)$ from Proposition \ref{prop:conj},
we obtain the claimed bound for
$\mathbf{E}(N(m,\alpha)) \leq m + 1 + T + T^*$.
\end{pf*}

We now turn to Theorem \ref{corol:1} on the expected value of
$\ln{\mathscr C}(A)$. In Theorem~\ref{th:main_lp} we derived tail
estimates on the probability distribution of the GCC condition number.
For the sake of clarity, we include the following simple observation
showing how to use these tail estimates to bound the expected value of
the logarithm of the condition number.
\begin{proposition}\label{prop:E}
Let $Z$ be a random variable, almost surely greater or equal than 1,
satisfying, for some $K,t_0>0$, that $\mathsf{Prob}\{Z\geq t\}\leq K
t^{-1}$ for all $t\ge t_0$. Then $\mathbf{E}(\ln Z)\leq\ln t_0
+\frac{K}{t_0}$.
\end{proposition}
\begin{pf}
We have $\mathsf{Prob}\{\ln Z\geq s\}\leq K e^{-s}$ for all $s>\ln
t_0$. Therefore,
\[
\mathbf{E}(\ln Z) = \int_0^{\infty}\mathsf{Prob}\{\ln Z\geq s\} \,ds
\leq\ln t_0 + \int_{\ln t_0}^{\infty} K e^{-s} \,dt = \ln t_0 +
\frac{K}{t_0}
\]
as claimed.
\end{pf}

We next proceed to prove Theorem \ref{corol:1}.
To simplify notation we put
\begin{eqnarray*}
P_{n,m}(\varepsilon) &:=&
\mathsf{Prob}\{A\in\mathcal{I}_{n,m}\mbox{ and }
{\mathscr C}(A)\ge\varepsilon^{-1}\},\\
Q_{n,m}(\varepsilon) &:=&
\mathsf{Prob}\{A\in\mathcal{F}_{n,m}\mbox{ and }
{\mathscr C}(A)\ge\varepsilon^{-1}\}.
\end{eqnarray*}
\begin{lemma}\label{lem:inf}
For any $n> m\geq1$ and $\varepsilon\in(0,1]$ we have:
\begin{longlist}
\item
If $\varepsilon^{-1}\ge13(m+1)^{3/2}$ then
$P_{n,m}(\varepsilon)\le2e (m+1)^{3/2} \varepsilon$.
\item
If $\varepsilon^{-1}\ge(m+1)^{2}$ then
$Q_{n,m}(\varepsilon)\le\sqrt{2\pi e} (m+1)^{7/4} \varepsilon$.
\end{longlist}
\end{lemma}
\begin{pf} (i)
Theorem \ref{th:main_lp} tells us that
\[
P_{n,m}(\varepsilon) \le\pmatrix{n\cr{m+1}} C(m,m)
\int_0^\varepsilon(1-t^2)^{({m^2-2})/{2}}\bigl(1-\lambda
_m(t)\bigr)^{n-m-1} \,dt.
\]
Recall formula (\ref{eq:volcap}) for the relative volume $\lambda_m(t)$
of a cap of radius $\arccos(t)$ on $S^m$. Recall also from Lemma
\ref{le:easy} that $\alpha_m := \frac{2\mathcal O_{m-1}}{\mathcal O_m}
\leq\sqrt{m}$. The first order derivative of $\lambda_m(t)$
\[
\frac{d\lambda_m(t)}{dt} = -\frac12\alpha_m (1-t^2)^{({m-2})/{2}}
\]
is increasing; hence $\lambda_m$ is a convex function. Moreover,
$\lambda_m(0)=1/2$. This implies $2\lambda_m(t)\ge1-\alpha_m t$ for all
$t\in[0,1]$; hence $1-\lambda_m(t)\le\frac12(1+\alpha_mt)$.

Bounding $C(m,m)$ as in Proposition \ref{prop:conj} we arrive at the estimate
%
%
\begin{equation}\label{eq:star}
P_{n,m}(\varepsilon) \le2(m+1)\sqrt{m} \frac1{2^{n}}{n\choose{m+1}}
\bigl(1+\sqrt{m}\varepsilon\bigr)^{n-m-1} \varepsilon.
\end{equation}
We now proceed dividing by cases. Suppose that
$\varepsilon^{-1}\ge13(m+1)^{3/2}$.
\begin{Case}[{[$n\le13 (m+1)$]}]\label{case1}
In this case $\varepsilon^{-1}\ge n\sqrt{m}$ and hence, using
(\ref{eq:star}),
\[
P_{n,m}(\varepsilon) \le2(m+1)\sqrt{m} (1+ 1/n)^{n-m-1} \varepsilon
\le2e(m+1)\sqrt{m} \varepsilon.
\]
\end{Case}
\begin{Case}[{[$n> 13(m+1)$]}]\label{case2}
This implies $\ln(e \frac{n}{m+1} ) \leq\frac{n}{m+1} \ln(\frac{4}{3}
)$, and it follows that
%
%
\begin{equation}\label{eq:P1}
\pmatrix{n\cr{m+1}} \le\frac1{(m+1)!} n^{m+1}
\le\biggl(\frac{en}{m+1} \biggr)^{m+1}\leq\biggl(\frac{4}{3} \biggr)^n.
\end{equation}
Since, in addition, $\varepsilon^{-1}\geq13(m+1)\sqrt{m}\ge2\sqrt{m}$
we get from (\ref{eq:star})
\begin{eqnarray*}
P_{n,m}(\varepsilon) \le2(m+1)\sqrt{m} \frac1{2^n}\pmatrix{n\cr{m+1}}
\biggl(\frac{3}{2} \biggr)^{n} \varepsilon
\le2(m+1)\sqrt{m} \varepsilon.\\[-15pt]
\end{eqnarray*}
\end{Case}

(ii)
Theorem \ref{th:main_lp} implies that
\begin{eqnarray*}
Q_{n,m}(\varepsilon) &=& \sum_{k=1}^m \pmatrix{n\cr{k+1}}C(m,k)
\int_0^\varepsilon t^{m-k} (1-t^2)^{km/2-1}
\lambda_m(t)^{n-k-1} \,dt\\
&\leq& \sum_{k=1}^m \pmatrix{n\cr{k+1}}C(m,k)
\varepsilon^{m-k+1} 2^{-(n-k-1)}\\
&\leq& \sum_{k=1}^m C(m,k)
\varepsilon^{m-k+1} 2^{k+1} ,
\end{eqnarray*}
the second line since $\lambda_m(t)\le\frac12$.
Using Proposition \ref{prop:conj} we obtain
\begin{eqnarray*}
Q_{n,m}(\varepsilon)
&\le& \varepsilon\sqrt{2\pi} (m+1)^{7/4} \sum_{k=1}^m \sqrt
{\pmatrix{m\cr
{k}}}
\bigl( (m+1)\varepsilon\bigr)^{m-k} \\
&\le& \varepsilon\sqrt{2\pi} (m+1)^{7/4} \sum_{k=1}^m \pmatrix
{{m}\cr{k}}
\bigl( (m+1)\varepsilon\bigr)^{m-k} \\
&\leq& \varepsilon\sqrt{2\pi} (m+1)^{7/4} \bigl(1+(m+1)\varepsilon\bigr)^m.
\end{eqnarray*}
Under the assumption $\varepsilon^{-1}\ge(m+1)^2$
we have $(m+1)\varepsilon\le\frac1{(m+1)}$, and hence
\[
Q_{n,m}(\varepsilon) \le\varepsilon\sqrt{2\pi} (m+1)^{7/4}\sqrt{e}.
\]
\upqed\end{pf}
\begin{pf*}{Proof of Theorem \protect\ref{corol:1}}
For $\varepsilon^{-1}\ge13(m+1)^2$ we have, by Lemma \ref{lem:inf},
\begin{eqnarray*}
\mathsf{Prob}\{{\mathscr C}(A) \geq\varepsilon^{-1} \}
&=& P_{n,m}(\varepsilon) + Q_{n,m}(\varepsilon)\\
&\leq& \bigl(2e(m+1)^{3/2}+\sqrt{2\pi e}(m+1)^{7/4}\bigr) \varepsilon\\
&\leq& 9.6 (m+1)^2 \varepsilon.
\end{eqnarray*}
An application of Proposition \ref{prop:E}
with $K = 9.6 (m+1)^2$ and
$t_0=13(m+1)^2$ shows that
\[
\mathbf{E}(\ln{\mathscr C}(A))\leq2\ln(m+1) +\ln13 +9.6/13 \leq
2\ln(m+1) +3.31.
\]
\upqed\end{pf*}

\subsection{On calculating the $C(m,k)$}\label{sec:Cmk_part}

We describe a method for calculating the $C(m,k)$.
For $1\le k\le m< n$ we define the following integrals:
\[
I(n,m,k):= 2^{n-1}\pmatrix{n\cr k+1}
\int_0^1 t^{m-k} (1-t^2)^{km/2-1}\lambda_m(t)^{n-k-1} \,dt.
\]
By setting $\varepsilon=1$
in the first part of Theorem \ref{th:main_lp} we
get from (\ref{eq:wendel}) that, for all $n> m$,
%
%
\begin{equation}\label{eq:lins}
\sum_{k=1}^m I(n,m,k)C(m,k)
= \sum_{k=0}^{m-1}\pmatrix{n-1\cr k}.
\end{equation}
By taking $m$ different values of $n$ (e.g., $n=m+1,\ldots,2m$)
one obtains a (square) system of linear equations in the $C(m,k)$.
Solving this system with Maple (symbolically for even $m$ and
numerically for odd $m$) we obtained Table \ref{table:1}.

We can further use (\ref{eq:lins}) to obtain expressions for
$C(m,k)$ for values of $k$ other than $1$ and $m$. We do so
for $k=m-1$.
\begin{proposition}
For all $m\geq2$,
\[
C(m,m-1)=\frac{m(m-1)}{2^{m-1}}(1+\alpha_m^2)
\qquad\mbox{where }
\alpha_m = \frac{2\mathcal O_{m-1}}{\mathcal O_m} .
\]
\end{proposition}
\begin{pf*}{Sketch of proof}
Put $J(n,m,k):=\int_0^1 t^{m-k} (1-t^2)^{km/2-1} \lambda _m(t)^{n-k-1}
\,dt$ so that $I(n,m,k)=2^{n-1}{n\choose k+1}J(n,m,k)$. In the
following we write $N=n-m$. One can prove that for fixed $m$ the
following asymptotic expansion holds for $N\to\infty$:
\begin{eqnarray*}
2^{n-m-1} J(n,m,m) &=& \frac1{\alpha_m}\frac1{N} -
\frac{m(m-1)}{\alpha_m^3}\frac1{N^3}
+ \mathcal O\biggl(\frac1{N^5} \biggr) ,\\
2^{n-m} J(n,m,m-1) &=& \frac1{\alpha_m^2} \frac1{(N+1)(N+2)} + \mathcal
O\biggl(\frac1{N^4} \biggr) .
\end{eqnarray*}
It follows after a short calculation that the left-hand side of
(\ref{eq:lins}) has the following expansion:
\begin{eqnarray*}
&&C(m,m)\frac{2^m}{(m+1)!}
\biggl(\frac1{\alpha_m} N^m + \frac{a_1(m)}{\alpha_m} N^{m-1}
+ \biggl(\frac{a_2(m)}{\alpha_m} - \frac{m(m-1)}{\alpha_m^3} \biggr)
N^{m-2} \biggr) \\
&&\qquad{} + C(m,m-1) \frac{2^{m-1}}{m!}\frac1{\alpha_m^2} N^{m-2}
+ \mathcal O(N^{m-3} ),
\end{eqnarray*}
where
\begin{eqnarray*}
a_1(m) &:=& \sum_{0\le j\le m} j = \frac12 m(m+1),\\
a_2(m) &:=& \sum_{0\le i<j\le m} ij = \frac1{24} (m+1)m(m-1)(m-2)(3m+2).
\end{eqnarray*}
Now we expand the right-hand side of (\ref{eq:lins}) to obtain
\begin{eqnarray*}
&&\frac1{m!} N^m + \biggl(\frac{a_1(m-1)}{m!} + \frac1{(m-1)!} \biggr) N^{m-1} \\
&&\qquad{} + \biggl(\frac{a_2(m-1)}{m!}
+ \frac{a_1(m-1)}{(m-1)!}+\frac1{(m-2)!} \biggr) N^{m-2}\\
&&\qquad{} + \mathcal O(N^{m-3}).
\end{eqnarray*}
By comparing the coefficients of $N^m$ (or those of $N^{m-1}$) we obtain
\[
C(m,m) = \frac{m+1}{2^m}\alpha_m .
\]
By comparing the coefficients of $N^{m-2}$ we get,
after a short calculation,
\begin{eqnarray*}
C(m,m-1) &=&
\frac{\mathcal O_{m-1}^2}{\mathcal O_m^2 2^{m-3}} \biggl(
a_2(m-1)-a_2(m)+ma_1(m-1) \\
&&\hspace*{69.5pt}{} + m(m-1) +
\frac{m(m-1)\mathcal O_m^2}{4\mathcal O_{m-1}^2} \biggr),
\end{eqnarray*}
and simplifying this expression, the claimed result
follows.
\end{pf*}

Finding a closed form for all coefficients $C(m,k)$
remains a challenging task.

\section*{Acknowledgment}
We thank Dennis Amelunxen for pointing out an error
in a previous version in the paper.

%

%
\printaddresses


\begin{thebibliography}{31}

\bibitem{Agmon}
%
\begin{barticle}[mr]
\bauthor{\bsnm{Agmon},~\bfnm{Shmuel}\binits{S.}}
(\byear{1954}).
\btitle{The relaxation method for linear inequalities}.
\bjournal{Canad. J. Math.}
\bvolume{6}
\bpages{382--392}.
\bid{mr={0062786}}
\end{barticle}
%
\endbibitem

\bibitem{AmBu08}
%
\begin{bmisc}[auto:SpringerTagBib|2009-01-14|16:51:27]
\bauthor{\bsnm{B\"{u}rgisser},~\bfnm{P.}\binits{P.}} \AND
\bauthor{\bsnm{Amelunxen},~\bfnm{D.}\binits{D.}}
(\byear{2008}).
\bhowpublished{Uniform smoothed analysis
of a condition number for linear programming}. Accepted for \textit{Math. Program. A.}
Available at
\href{http://www.arxiv.org/abs/0803.0925}{arXiv:0803.0925}.
\end{bmisc}
%
\endbibitem

\bibitem{ChC00}
%
\begin{barticle}[mr]
\bauthor{\bsnm{Cheung},~\bfnm{Dennis}\binits{D.}} \AND
\bauthor{\bsnm{Cucker},~\bfnm{Felipe}\binits{F.}}
(\byear{2001}).
\btitle{A new condition number for linear programming}.
\bjournal{Math. Program.}
\bvolume{91}
\bpages{163--174}.
\bid{mr={1865268}}
\end{barticle}
%
\endbibitem

\bibitem{ChC01}
%
\begin{barticle}[mr]
\bauthor{\bsnm{Cheung},~\bfnm{D.}\binits{D.}} \AND
\bauthor{\bsnm{Cucker},~\bfnm{F.}\binits{F.}}
(\byear{2002}).
\btitle{Probabilistic analysis of condition numbers for linear programming}.
\bjournal{J. Optim. Theory Appl.}
\bvolume{114}
\bpages{55--67}.
\bid{doi={10.1023/A:1015460004163}, mr={1910854}}
\end{barticle}
%
\endbibitem

\bibitem{ChCH05}
%
\begin{barticle}[mr]
\bauthor{\bsnm{Cheung},~\bfnm{Dennis}\binits{D.}},
\bauthor{\bsnm{Cucker},~\bfnm{Felipe}\binits{F.}} \AND
\bauthor{\bsnm{Hauser},~\bfnm{Raphael}\binits{R.}}
(\byear{2005}).
\btitle{Tail decay and moment estimates of a condition number for
random linear
conic systems}.
\bjournal{SIAM J. Optim.}
\bvolume{15}
\bpages{1237--1261}.
\bid{doi={10.1137/S105262340343470X}, mr={2178497}}
\end{barticle}
%
\endbibitem

\bibitem{CP01}
%
\begin{barticle}[mr]
\bauthor{\bsnm{Cucker},~\bfnm{Felipe}\binits{F.}} \AND
\bauthor{\bsnm{Pe{\~n}a},~\bfnm{Javier}\binits{J.}}
(\byear{2002}).
\btitle{A primal-dual algorithm for solving polyhedral conic systems
with a
finite-precision machine}.
\bjournal{SIAM J. Optim.}
\bvolume{12}
\bpages{522--554}.
\bid{doi={10.1137/S1052623401386794}, mr={1885574}}
\end{barticle}
%
\endbibitem

\bibitem{CW01}
%
\begin{barticle}[mr]
\bauthor{\bsnm{Cucker},~\bfnm{Felipe}\binits{F.}} \AND
\bauthor{\bsnm{Wschebor},~\bfnm{Mario}\binits{M.}}
(\byear{2002}).
\btitle{On the expected condition number of linear programming problems}.
\bjournal{Numer. Math.}
\bvolume{94}
\bpages{419--478}.
\bid{doi={10.1007/s00211-002-0385-1}, mr={1981163}}
\end{barticle}
%
\endbibitem

\bibitem{DST}
%
\begin{bmisc}[auto:SpringerTagBib|2009-01-14|16:51:27]
\bauthor{\bsnm{Dunagan},~\bfnm{J.}\binits{J.}},
\bauthor{\bsnm{Spielman},~\bfnm{D.~A.}\binits{D.~A.}} \AND
\bauthor{\bsnm{Teng},~\bfnm{S.-H.}\binits{S.-H.}}
(\byear{2009}).
\bhowpublished{Smoothed analysis of condition numbers and complexity implications
for linear programming. \textit{Math. Programming.} To appear. Available at
\url{http://arxiv.org/abs/cs/0302011v2}.}
\end{bmisc}
%
\endbibitem

\bibitem{Dvor56}
%
\begin{barticle}[mr]
\bauthor{\bsnm{Dvoretzky},~\bfnm{Aryeh}\binits{A.}}
(\byear{1956}).
\btitle{On covering a circle by randomly placed arcs}.
\bjournal{Proc. Natl. Acad. Sci. USA}
\bvolume{42}
\bpages{199--203}.
\bid{mr={0079365}}
\end{barticle}
%
\endbibitem

\bibitem{Gilbert65}
%
\begin{barticle}[mr]
\bauthor{\bsnm{Gilbert},~\bfnm{E.~N.}\binits{E.~N.}}
(\byear{1965}).
\btitle{The probability of covering a sphere with {$N$} circular caps}.
\bjournal{Biometrika}
\bvolume{52}
\bpages{323--330}.
\bid{mr={0207005}}
\end{barticle}
%
\endbibitem

\bibitem{Goff80}
%
\begin{barticle}[mr]
\bauthor{\bsnm{Goffin},~\bfnm{J.-L.}\binits{J.-L.}}
(\byear{1980}).
\btitle{The relaxation method for solving systems of linear inequalities}.
\bjournal{Math. Oper. Res.}
\bvolume{5}
\bpages{388--414}.
\bid{doi={10.1287/moor.5.3.388}, mr={594854}}
\end{barticle}
%
\endbibitem

\bibitem{Hall85}
%
\begin{barticle}[mr]
\bauthor{\bsnm{Hall},~\bfnm{Peter}\binits{P.}}
(\byear{1985}).
\btitle{On the coverage of {$k$}-dimensional space by {$k$}-dimensional
spheres}.
\bjournal{Ann. Probab.}
\bvolume{13}
\bpages{991--1002}.
\bid{mr={799434}}
\end{barticle}
%
\endbibitem

\bibitem{Hall88}
%
\begin{bbook}[mr]
\bauthor{\bsnm{Hall},~\bfnm{Peter}\binits{P.}}
(\byear{1988}).
\btitle{Introduction to the Theory of Coverage Processes}.
\bpublisher{Wiley}, \baddress{New York}.
\bid{mr={973404}}
\end{bbook}
%
\endbibitem

\bibitem{HM06}
%
\begin{barticle}[mr]
\bauthor{\bsnm{Hauser},~\bfnm{Raphael}\binits{R.}} \AND
\bauthor{\bsnm{M{\"u}ller},~\bfnm{Tobias}\binits{T.}}
(\byear{2009}).
\btitle{Conditioning of random conic systems under a general family of input
distributions}.
\bjournal{Found. Comput. Math.}
\bvolume{9}
\bpages{335--358}.
\bid{doi={10.1007/s10208-008-9034-0}, mr={2496555}}
\end{barticle}
%
\endbibitem

\bibitem{Janson86}
%
\begin{barticle}[mr]
\bauthor{\bsnm{Janson},~\bfnm{Svante}\binits{S.}}
(\byear{1986}).
\btitle{Random coverings in several dimensions}.
\bjournal{Acta Math.}
\bvolume{156}
\bpages{83--118}.
\bid{doi={10.1007/BF02399201}, mr={822331}}
\end{barticle}
%
\endbibitem

\bibitem{Kahane59}
%
\begin{barticle}[mr]
\bauthor{\bsnm{Kahane},~\bfnm{Jean-Pierre}\binits{J.-P.}}
(\byear{1959}).
\btitle{Sur le recouvrement d'un cercle par des arcs dispos\'es au hasard}.
\bjournal{C. R. Math. Acad. Sci. Paris}
\bvolume{248}
\bpages{184--186}.
\bid{mr={0103533}}
\end{barticle}
%
\endbibitem

\bibitem{Miles68a}
%
\begin{barticle}[auto:SpringerTagBib|2009-01-14|16:51:27]
\bauthor{\bsnm{Miles},~\bfnm{R.~E.}\binits{R.~E.}}
(\byear{1968}).
\btitle{Random caps on a sphere}.
\bjournal{Ann. Math. Statist.}
\bvolume{39}
\bpages{1371}.
\end{barticle}
%
\endbibitem

\bibitem{Miles69}
%
\begin{barticle}[mr]
\bauthor{\bsnm{Miles},~\bfnm{R.~E.}\binits{R.~E.}}
(\byear{1969}).
\btitle{The asymptotic values of certain coverage probabilities}.
\bjournal{Biometrika}
\bvolume{56}
\bpages{661--680}.
\bid{mr={0254953}}
\end{barticle}
%
\endbibitem

\bibitem{mile71}
%
\begin{barticle}[mr]
\bauthor{\bsnm{Miles},~\bfnm{R.~E.}\binits{R.~E.}}
(\byear{1971}).
\btitle{Isotropic random simplices}.
\bjournal{Adv. in Appl. Probab.}
\bvolume{3}
\bpages{353--382}.
\bid{mr={0309164}}
\end{barticle}
%
\endbibitem

\bibitem{MoFa62}
%
\begin{barticle}[mr]
\bauthor{\bsnm{Moran},~\bfnm{P.~A.~P.}\binits{P.~A.~P.}} \AND
\bauthor{\bparticle{Fazekas~de }\bsnm{St.~Groth},~\bfnm{S.}\binits{S.}}
(\byear{1962}).
\btitle{Random circles on a sphere}.
\bjournal{Biometrika}
\bvolume{49}
\bpages{389--396}.
\bid{mr={0156434}}
\end{barticle}
%
\endbibitem

\bibitem{Motzkin}
%
\begin{barticle}[mr]
\bauthor{\bsnm{Motzkin},~\bfnm{T.~S.}\binits{T.~S.}} \AND
\bauthor{\bsnm{Schoenberg},~\bfnm{I.~J.}\binits{I.~J.}}
(\byear{1954}).
\btitle{The relaxation method for linear inequalities}.
\bjournal{Canad. J. Math.}
\bvolume{6}
\bpages{393--404}.
\bid{mr={0062787}}
\end{barticle}
%
\endbibitem

\bibitem{reitz}
%
\begin{barticle}[mr]
\bauthor{\bsnm{Reitzner},~\bfnm{Matthias}\binits{M.}}
(\byear{2002}).
\btitle{Random points on the boundary of smooth convex bodies}.
\bjournal{Trans. Amer. Math. Soc.}
\bvolume{354}
\bpages{2243--2278}.
\bid{doi={10.1090/S0002-9947-02-02962-8}, mr={1885651}}
\end{barticle}
%
\endbibitem

\bibitem{Rosenblatt}
%
\begin{bbook}[mr]
\bauthor{\bsnm{Rosenblatt},~\bfnm{Frank}\binits{F.}}
(\byear{1962}).
\btitle{Principles of Neurodynamics. {P}erceptrons and the Theory of Brain
Mechanisms}.
\bpublisher{Spartan Books}, \baddress{Washington, DC}.
\bid{mr={0135635}}
\end{bbook}
%
\endbibitem

\bibitem{sant76}
%
\begin{bbook}[mr]
\bauthor{\bsnm{Santal{\'o}},~\bfnm{Luis~A.}\binits{L.~A.}}
(\byear{1976}).
\btitle{Integral Geometry and Geometric Probability}.
\bpublisher{Addison-Wesley}, \baddress{Reading, MA}.
\bid{mr={0433364}}
\end{bbook}
%
\endbibitem

\bibitem{Siegel79}
%
\begin{barticle}[mr]
\bauthor{\bsnm{Siegel},~\bfnm{Andrew~F.}\binits{A.~F.}}
(\byear{1979}).
\btitle{Asymptotic coverage distributions on the circle}.
\bjournal{Ann. Probab.}
\bvolume{7}
\bpages{651--661}.
\bid{mr={537212}}
\end{barticle}
%
\endbibitem

\bibitem{SieHo82}
%
\begin{barticle}[mr]
\bauthor{\bsnm{Siegel},~\bfnm{Andrew~F.}\binits{A.~F.}} \AND
\bauthor{\bsnm{Holst},~\bfnm{Lars}\binits{L.}}
(\byear{1982}).
\btitle{Covering the circle with random arcs of random sizes}.
\bjournal{J. Appl. Probab.}
\bvolume{19}
\bpages{373--381}.
\bid{mr={649974}}
\end{barticle}
%
\endbibitem

\bibitem{solo78}
%
\begin{bbook}[mr]
\bauthor{\bsnm{Solomon},~\bfnm{Herbert}\binits{H.}}
(\byear{1978}).
\btitle{Geometric Probability}.
\bpublisher{SIAM},
\baddress{Philadelphia, PA}.
\bid{mr={0488215}}
\end{bbook}
%
\endbibitem

\bibitem{Stev39}
%
\begin{barticle}[mr]
\bauthor{\bsnm{Stevens},~\bfnm{W.~L.}\binits{W.~L.}}
(\byear{1939}).
\btitle{Solution to a geometrical problem in probability}.
\bjournal{Ann. Eugenics}
\bvolume{9}
\bpages{315--320}.
\bid{mr={0001479}}
\end{barticle}
%
\endbibitem

\bibitem{Wendel62}
%
\begin{barticle}[mr]
\bauthor{\bsnm{Wendel},~\bfnm{J.~G.}\binits{J.~G.}}
(\byear{1962}).
\btitle{A problem in geometric probability}.
\bjournal{Math. Scand.}
\bvolume{11}
\bpages{109--111}.
\bid{mr={0146858}}
\end{barticle}
%
\endbibitem

\bibitem{Whit1897}
%
\begin{bbook}[auto:SpringerTagBib|2009-01-14|16:51:27]
\bauthor{\bsnm{Whitworth},~\bfnm{W.~A.}\binits{W.~A.}}
(\byear{1965}).
\btitle{DCC Exercises in Choice and Chance}.
\bpublisher{Dover},
\baddress{New York}.
\end{bbook}
%
\endbibitem

\bibitem{zaehle}
%
\begin{barticle}[mr]
\bauthor{\bsnm{Z{\"a}hle},~\bfnm{M.}\binits{M.}}
(\byear{1990}).
\btitle{A kinematic formula and moment measures of random sets}.
\bjournal{Math. Nachr.}
\bvolume{149}
\bpages{325--340}.
\bid{mr={1124814}}
\end{barticle}
%
\endbibitem

\end{thebibliography}
\end{document}